\documentclass[final,3p,times]{elsarticle}

\usepackage{amssymb}
\usepackage{amsthm}
\usepackage{amsmath}
\usepackage{amsfonts}

\newtheorem{thm}{Theorem}

\newtheorem{lem}{Lemma}

\newtheorem{defn}{Definition}

\newtheorem*{dem}{Proof}

\begin{document}

\begin{frontmatter}



\title{Some Results on the $1$-Laplacian Elliptic Problems with Singularities and Robin Boundary Conditions}

\author[label1]{Mohamed El Hichami}
\ead{elhichami.7@gmail.com}
\address[label1]{{Laboratory LIPIM, National School of Applied Sciences Khouribga},
            { Sultan Moulay Slimane University},
            {Khouribga},
            {25000},
            {Morocco}}

\author[label1]{Youssef El Hadfi}
\ead{yelhadfi@gmail.com}

\begin{abstract}
In this paper, we investigate the existence and uniqueness of solutions for the following model problem, involving singularities and inhomogeneous Robin boundary conditions
\begin{equation*}
    \left\{
  \begin{array}{ll}
 -\Delta_{p}u_{p}=\frac{f}{u_{p}^{\gamma}}& \hbox{in $\Omega,$} \\
  \frac{\partial u_{p}}{\partial \sigma}+\lambda\vert u_{p}\vert^{p-2} u_{p}+\vert u_{p}\vert^{s-1}u_{p}=\frac{g}{u_{p}^{\eta}} & \hbox{on $\partial\Omega,$} \\
  \end{array}
\right.
\end{equation*}
where $\Omega \subset \mathbb{R}^{m}$ represents an open bounded domain, with smooth boundary, $m \geq 2$, the symbol  $\sigma $ stands for the unit outward normal vector, $ \Delta_{p}u:=\mbox{div}(\vert\nabla u\vert^{p-2}\nabla u) $ is the $p-$Laplacian operator $(1\leq p<m),$ consider  $0<\gamma\leq 1,$ $ \eta>0$ and  $s\geq 1.$ The function  $ f\in L^{\frac{m}{p}}(\Omega)$ is a nonnegative additionally $ \lambda$ and $ g$ are nonnegative functions in $L^{\infty}(\partial \Omega).$
\end{abstract}
\begin{keyword}
 $1-$Laplacian operator\sep Singular term \sep elliptic equations \sep Robin boundary conditions \sep Functions of bounded variations.
\end{keyword}

\end{frontmatter}


\section{Introduction}
\par The purpose of this paper is to examine the following issue
\begin{equation}\label{53}
    \left\{
  \begin{array}{ll}
   -\Delta_{p}u_{p}=\frac{f}{u_{p}^{\gamma}}& \hbox{in $\Omega,$} \\
  \frac{\partial u_{p}}{\partial \sigma}+\lambda\vert u_{p}\vert^{p-2} u_{p}+\vert u_{p}\vert^{s-1}u_{p}=\frac{g}{u_{p}^{\eta}} & \hbox{on $\partial\Omega,$} \\
  \end{array}
\right.
\end{equation}
where $1\leq p<m,$ $0<\gamma\leq 1,$ $ \eta>0,$  $s\geq 1,$ $\Omega \subset \mathbb{R}^{m (m \geq 2)}$ represents an open bounded domain with smooth boundary, the symbol  $\sigma $ stands for the unit outward normal vector, $ 0 \leq f\in L^{\frac{m}{p}}(\Omega)$ and $ \lambda,$  $ g$ are nonnegative functions in $L^{\infty}(\partial \Omega).$
\par The 1-Laplacian operator, denoted as
\[
\Delta_1 u := \operatorname{div}\left(\frac{D u}{|D u|}\right),
\]
plays a crucial role in various mathematical theories and practical applications, especially in fields where preserving sharp edges and fine details is essential. It is prominently used in image processing for tasks like denoising and segmentation. The 1-Laplacian is the core operator in the celebrated Total Variation (TV) denoising model introduced by Rudin, Osher, and Fatemi \cite{ROF}, which is widely applied in medical imaging, such as MRI (see \cite{BASU, TRI}), to preserve edges while reducing noise. Unlike the p-Laplacian operator, which smooths fine details, the $1$-Laplacian enables isotropic diffusion within level surfaces without crossing them, effectively preserving edges.

Beyond medical imaging, the 1-Laplacian operator has been applied in other domains, such as synthetic aperture radar (SAR) imaging, to address speckle noise \cite{Lee}. Variational formulations involving the 1-Laplacian operator provide faithful solutions that retain key features of the original image, making it ideal for edge-preserving techniques. It is also useful in solving inverse ill-posed problems, where the goal is to recover meaningful information from incomplete or noisy data.

The mathematical properties of the $1$-Laplacian have also been studied in the context of nonlinear elliptic and parabolic equations, particularly as a limit case of the p-Laplacian for \( p=1 \). Its ability to handle discontinuous solutions in the space of functions of bounded variation makes it suitable for advanced applications, such as object recognition, super-resolution, and image classification.

Physical applications manifest in addressing problems such as \eqref{53}, with illustrative instances found in signal transmissions \cite{60222} and the theory of non-Newtonian pseudo-plastic fluids \cite{8126}. Furthermore, additional motivations emerge in fields like image restoration \cite{6012} and game theory \cite{7022}. Equation \eqref{53} can be regarded as a broader version of the classical p-Laplacian equation
\begin{equation}\label{52}
    \left\{
  \begin{array}{ll}
   -\Delta_{p}u_{p}=fh_{1}(u_{p})& \hbox{in $\Omega,$} \\
  \frac{\partial u_{p}}{\partial \sigma}+\lambda\vert u_{p}\vert^{p-2} u_{p}+k(u_{p})=gh_{2}(u_{p}) & \hbox{on $\partial\Omega,$} \\
  \end{array}
\right.
\end{equation}
with $ h_{1},$ $ h_{2}$ are continuous functions on $[0,+\infty[$ (possibly unbounded near zero) and $ k(t)=\vert t \vert^{s-1}t,$ $ s\geq 1.$
\par The problem \eqref{52}, specifically for the condition  $p>1,$ has been extensively researched. In \cite{6010},  when  $1<p<m,$ $h_{1}(t)=\frac{1}{t^{\gamma}},$ $0<\gamma<1,$ $ f\in L^{\zeta}(\Omega),$  $ \zeta=\frac{mp}{m(p-1)+p+\gamma(m-p)}$  and the Dirichlet condition $u=0$ replaces the second equation's Robin boundary conditions on $\partial \Omega,$ the authors showed that problem \eqref{52} has a weak solution (see also \cite{6014}).  Moreover in \cite{8133} it is considered the problem \eqref{52} for  $ p>1,$ $ h_{1}(t)=h_{2}(t)=1,$ $ k(t)=0,$ $f$ belongs to $L^{m,\infty}(\Omega),$ \(\lambda\) is a non-negative function in \(L^{\infty}(\partial \Omega)\) that is not identically zero, and $g$ is an element of \(L^{\infty}(\partial \Omega)\), it is proven that there is a weak solution. In the case $ p=2,$ $ h_{1}(t)=1,$ $ h_{2}(t)=\frac{1}{t^{\eta}},$ $\eta>0,$ $ k(t)=0,$ $ 0\leq f\in L^{\frac{2m}{m+2}}(\Omega),$ $0\leq \lambda \in L^{\infty}(\partial \Omega),$ $0 \leq g\in L^{r}(\partial \Omega)$ where $ r=\max\bigg(\frac{2(m-1)}{m+\eta(m-2)},1\bigg),$ problem \eqref{52} has been analyzed in \cite{8134} with the authors show the existence of a weak solution. The $p-$Laplacian operator was studied by several authors see \cite{6010,8132,6023,6024,Sbai1,Sbai2,Sbai3,Sbai4} and references therein.
\par When $ p=1,$ we denote by $ \Delta_{1}u:=\textnormal{div }\bigg(\frac{Du}{\vert Du\vert}\bigg)$ the $1-$Laplacian operator. The appropriate energy space for dealing with problems involving $1-$Laplacian is $ BV$ the space of functions of bounded variation. The notion of a proper solution is introduced in \cite{6003, 6004}.  Let us cite some works discussing the problems as \eqref{52}. If $ p=1,$ $ f\in L^{m}(\Omega),$  $h_{1}(t)=\frac{1}{t^{\gamma}},$ $0<\gamma\leq 1$ and
$ u=0$ on $ \partial \Omega,$ the authors in \cite{6017} showed the existence of a bounded solution to problem \eqref{52}. Furthermore in \cite{7028} it has been established that the existence, nonexistence and uniqueness  of solutions to \eqref{52} when $p=1,$ $f\in L^{1}(\Omega),$ $ u=0$ in $ \partial \Omega$ with $h_{1}:[0;+\infty[\rightarrow ]0;+\infty]$ is continuous, finite outside the origin, $ h_{1}(0)\neq 0,$
\begin{equation*}
  \exists \, \zeta,\hspace*{0.2cm} \gamma,\hspace*{0.2cm} \delta>0 \hspace*{0.2cm} \textnormal{such} \hspace*{0.2cm}\textnormal{that} \hspace*{0.2cm}h_{1}(t)\leq \frac{\zeta}{t^{\gamma}}\hspace*{0.2cm}\textnormal{if}\hspace*{0.2cm} t\leq \delta,
\end{equation*}
and $\underset{t\rightarrow\infty}{\lim}h_{1}(t):=h_{1}(\infty)<\infty.$  In the case of
$ p=1,$ $ f\in L^{m,\infty}(\Omega),$ $ g\in L^{\infty}(\partial \Omega),$ $\lambda \in L^{\infty}(\partial \Omega)$ and
$$M(f,g,\lambda)=\underset{u\in W^{1,1}(\Omega)\setminus \{0\}}{\sup}\frac{\int_{\Omega}fudy+\int_{\partial \Omega}gud\mathcal{H}^{m-1}}{\Vert u \Vert_{\lambda,1}}\leq 1,$$ the authors in \cite{8133} established the existence of a solution to  \eqref{52}. We refer to \cite{6016,6005,8113,8114,6019,6020} for more information concerning $1-$Laplacian operator.
\par When $1<p<m,$ the major difficulties in this study stem from the simultaneous presence of two distinct singularities, $\frac{1}{u^{\gamma}_{p}}$ and $\frac{1}{u^{\eta}_{p}},$ extending respectively over the domain $\Omega$ and the boundary $\partial \Omega,$ particularly on the set $\{u_{p}=0\}.$ To address this issue, we employ an approximation to regularize \eqref{1} and leverage the J. Leray and J. L. Lions Theorem. Additionally, the challenges are accentuated by the complexities associated with the Robin boundary conditions, which introduce a further layer of intricacy to the problem. As we approach the limit as $n\rightarrow +\infty$ in equation \eqref{9}, it is crucial to demonstrate the existence of a constant $\beta_{p}>0$ such that $ u_{n}>\beta_{p}$ on $ \partial \Omega.$ In the specific case of $ p=1,$ a notable difficulty arises when taking $p\rightarrow 1^{+}$ in equation \eqref{55}, as conventional methods cannot be applied to unbounded $BV-$functions. Therefore, we justify the use of truncated $BV-$functions.
\par This paper is organized as follows. Section 2 introduces some preliminaries and technical results. In Section 3, we prove the existence of solutions  for the problem where the Leray-Lions operator is represented by the p-Laplacian operator in the case $p>1.$
In Section 4, for $0<f\in L^{m}(\Omega)$ and $p=1,$ We investigate the existence and uniqueness of solutions. Finally, in Section \ref{sec5}, under the same assumptions, with the only difference being that the function $f$ is allowed to attain zero values, we show the existence of solutions.
\section{Preliminaries and technical results}
\par Let $\Omega$ be a bounded open subset of $\mathbb{R}^{m (m\geq 2)},$ with a smooth boundary. Given a function $ \kappa,$ we define  $\kappa^{+}:=\max(\kappa,0)$ and  $\kappa^{-}:=-\min(\kappa,0).$ Let \( \mathfrak{E} \subset \mathbb{R}^{m} \). We use \( \chi_{\mathfrak{E}} \) to represent the characteristic function of \( \mathfrak{E} \), and \( |\mathfrak{E}| \) as its Lebesgue measure.
Furthermore  $L^{\xi}(\mathfrak{E})$ is  the Lebesgue space of $\xi$-integrable functions on $\mathfrak{E}.$ Use $\mathcal{H}^{m-1}(\mathfrak{E})$ to indicate the $(m-1)$-dimensional Hausdorff measure.
We will utilize the standard Sobolev space $ W^{1,p}(\Omega)$ and $ \mathcal{M}(\Omega),$ which  represents the space of Radon measures with finite total variation throughout $ \Omega.$ Consider $ BV(\Omega)$ as the space of functions possessing bounded variation in $ \Omega.$  $u^{\ast}$ signifies the precise representative of $ u.$  For every $\tau>0$ and $r\in \mathbb{R},$ we define the auxiliary functions as usual
\begin{equation*}
    T_{\tau}(r):=\max(-\tau,\min(r,\tau)),\hspace*{0.1cm}
 \mathfrak{S}_{\tau}(r):=   \left\{
      \begin{array}{ll}
        1 & \hbox{$r\leq \tau$,} \\
        \frac{2\tau-r}{\tau} & \hbox{$\tau<r<2\tau$,} \\
       0  & \hbox{$ 2\tau\leq r,$}
      \end{array}
    \right.\mbox{ and } G_{\tau}(r):=1-\mathfrak{S}_{\tau}(r).
\end{equation*}
\par Given $ \lambda\in L^{\infty}(\partial \Omega)$ with $ \lambda\geq0$ and $1\leq p\leq m,$ we employ in $W^{1,p}(\Omega)$ the norm that is described as
\begin{equation}\label{47}
    \Vert u\Vert^{p}_{\lambda,p}=\int_{\Omega}\vert\nabla u \vert^{p}dy+\int_{\partial \Omega}\lambda\vert u\vert^{p} d\mathcal{H}^{m-1},\hspace*{0.5cm}u\in W^{1,p}(\Omega).
\end{equation}
In \cite{6022} proved that \eqref{47} is equivalent to the usual norm of $ W^{1,p}(\Omega).$
\par By \cite[Theorem 4.2]{6022} we underline that there exists constant $ C_{1}>0$ where
\begin{equation}\label{46}
\Vert u\Vert_{L^{\frac{(m-1)p}{m-p}}(\partial \Omega)}\leq C_{1}\Vert u \Vert_{W^{1,p}(\Omega)},  \hspace*{0.4cm}\mbox{ for all } u\in W^{1,p}(\Omega).
\end{equation}
Moreover if $\varrho<\frac{(m-1)p}{m-p} $ the mapping $W^{1,p}(\Omega)\rightarrow L^{\varrho}(\partial \Omega) $ is compact (see \cite[Theorem 6.2]{6022}).
We signify by $\mathcal{DM}^{\infty}(\Omega)$ the set of vector fields $ z$ belongs to $L^{\infty}(\Omega; \mathbb{R}^{m})$ and satisfy the condition that their
$\textnormal{div }z$ is an element of $\mathcal{M}(\Omega).$
$\mathcal{DM}_{loc}^{\infty}(\Omega)$ stands for the vector fields $z\in L^{\infty}(\Omega;\mathbb{R}^{m})$ where $\textnormal{div }z$ is a member of $\mathcal{M}(\omega)$ for all $ \omega \subset\subset \Omega.$
\par Let a  function $ u\in BV(\Omega)$ that means $ u\in L^{1}(\Omega)$ and whose gradient $Du$ is a vector  Radon measure  characterized by finite variation.
The space \( BV(\Omega) \) becomes a Banach space once it is equipped with one of the following two norms:
\begin{equation*}
   \Vert u\Vert_{BV(\Omega)}:=\int_{\Omega}\vert Du \vert +\int_{\Omega}\vert u\vert dy,
\mbox{ or }  \Vert u\Vert_{BV(\Omega)}:= \int_{\Omega}\vert Du \vert+\int_{\partial \Omega}\vert u\vert d\mathcal{H}^{m-1}.
\end{equation*}
 We highlight that the functional defined as
\begin{equation*}
u\mapsto \int_{\Omega}\varphi\vert Du \vert \hspace*{0.3cm}\textnormal{with}\hspace*{0.3cm} 0\leq\varphi \in C^{1}_{0}(\Omega),
\mbox{ and  }
u\mapsto \int_{\Omega}\vert Du \vert+\int_{\partial \Omega }\vert u\vert d\mathcal{H}^{m-1},
\end{equation*}
are lower semicontinuity in $ L^{1}(\Omega).$  For more detailed information concerning the $ BV$ function we refer to \cite{6002, 7042}.
\par Anzellotti's in \cite{6006} introduced  the following distribution $(z,Du):C^{\infty}_{0}(\Omega)\rightarrow \mathbb{R} $ as
\begin{equation}\label{95}
    <(z,Du),\varphi>:=-\int_{\Omega}u^{*}\varphi \mbox{div }z-\int_{\Omega}uz\nabla \varphi dy, \quad \varphi \in C^{\infty}_{0}(\Omega).
\end{equation}
with $\mbox{div }z\in L^{1}(\Omega)$ and $ u\in BV(\Omega)\cap L^{\infty}(\Omega)$    or  $z\in \mathcal{DM}^{\infty}(\Omega)$ and $ u\in BV(\Omega)\cap L^{\infty}(\Omega)\cap C(\Omega).$  Moreover the  authors in \cite{7034} extended the definition of \eqref{95} concerning the case in which $z\in \mathcal{DM}^{\infty}(\Omega)$ and $ u\in BV(\Omega)\cap L^{\infty}(\Omega).$ Finally the following  Lemmas are shown in \cite{6017}.
\begin{lem}
 If $u$ is an element of  $BV_{loc}(\Omega)\cap L^{\infty}(\Omega)$ and $ z$ belongs to  $\in \mathcal{DM}_{loc}^{\infty}(\Omega).$ Hence the functional $ (z,Du)$ is a Radon measure in $\Omega$ satisfying
\begin{equation*}
    \vert <(z,Du),\Phi> \vert\leq \Vert \Phi \Vert_{L^{\infty}(\Omega)}\Vert z\Vert_{L^{\infty}(\Upsilon)}\int_{\Upsilon}\vert Du\vert,
\end{equation*}
for every open set $ \Upsilon\subset\subset\Omega$ and for all $\Phi\in C^{\infty}_{0}(\Omega).$
\end{lem}
\begin{lem}\label{106}
Assume that $z$ belongs to $\mathcal{DM}_{\text{loc}}^{\infty}(\Omega)$ and $u$ belongs to $BV_{\text{loc}}(\Omega) \cap L^{\infty}(\Omega),$
 with $u^{\ast}\in L^{1}(\Omega, divz).$ Hence the following Green's Formula holds
\begin{equation*}
   \int_{\Omega}u^{\ast}\textnormal{div }z +\int_{\Omega}(z,Du)=\int_{\partial\Omega}[uz,\sigma]d\mathcal{H}^{m-1}.
\end{equation*}
\end{lem}

\section{Weak solution for $0\leq f\in L^{\frac{m}{p}}(\Omega)$ and $p>1$}

In this section, we investigate the existence of a weak solution for the problem presented in \eqref{53} involving the Leray-Lions operators.

Let \( \Omega \) be an open, bounded subset of \( \mathbb{R}^{m} \) with \( m \geq 2 \) and a smooth boundary. Let  $ \aleph:\Omega\times \mathbb{R}^{m}\rightarrow \mathbb{R}^{m} $ be a Carath\'{e}odory function. For every $y\in\Omega$ and $\forall \hbar, \hslash \in \mathbb{R}^{m}$ such that
\begin{align}
   & \aleph(y,\hbar).\hbar \geq \alpha \vert\hbar\vert^{p}, \hspace{0.5cm} \alpha >0,\label{b}\\
   &\vert \aleph(y,\hbar)\vert\leq l(y)+\varpi \vert\hbar\vert^{p-1},\hspace{0.5cm}  \varpi>0,\label{c} \\
   & (\aleph(y,\hbar)-\aleph(y,\hslash)).(\hbar-\hslash)>0, \hspace{0.5cm}\hbar \neq \hslash, \label{d}
\end{align}
where $1<p<m$ and $l\in L^{p^{\prime}}(\Omega),  p^{\prime}=\frac{p}{p-1}.$
\par Let $1<p<m,$ we consider the nonlinear problem
\begin{equation}\label{1}
    \left\{
  \begin{array}{ll}
   -\mbox{div }(\aleph(y, \nabla u_{p}))=\frac{f}{u_{p}^{\gamma}}& \hbox{in $\Omega,$} \\
  \frac{\partial u_{p}}{\partial \sigma}+\lambda\vert u_{p}\vert^{p-2} u_{p}+\vert u_{p}\vert^{s-1}u_{p}=\frac{g}{u_{p}^{\eta}} & \hbox{on $\partial\Omega,$} \\
  \end{array}
\right.
\end{equation}
with $ f\in L^{\frac{m}{p}}(\Omega)$ is a nonnegative function and $ \lambda,$  $ g$ are nonnegative functions in $L^{\infty}(\partial \Omega).$

\par We present the solution concept for this case.
\begin{defn}
A function $u_{p}\in W^{1,p}(\Omega)$ is a weak solution to \eqref{1}, if $ \frac{g}{u^{\eta}_{p}}\in L^{1}(\partial \Omega),$ $ \vert u_{p}\vert^{s-1}u_{p}\in L^{1}(\partial \Omega),$ $ \frac{f}{u_{p}^{\gamma}}\in L^{1}(\Omega)$ and it satisfies
\begin{equation}\label{19}
\int_{\Omega}\aleph(x,\nabla u_{p})\nabla \varphi dy+\int_{\partial \Omega} \lambda \vert u_{p}\vert^{p-2}u_{p}\varphi d\mathcal{H}^{m-1}+\int_{\partial \Omega}\vert u_{p}\vert^{s-1}u_{p}\varphi d\mathcal{H}^{m-1}=\int_{\Omega}\frac{f}{u^{\gamma}_{p}}\varphi dy+\int_{\partial \Omega}\frac{g}{u_{p}^{\eta}}\varphi d\mathcal{H}^{m-1},
\end{equation}
for all $\varphi \in W^{1,p}(\Omega)\cap L^{\infty}(\Omega).$
\end{defn}
This allows us to provide our first result.
\begin{thm}\label{48}
  Consider $ 0\leq f\in L^{\frac{m}{p}}(\Omega),$ let  $\lambda$ and $ g$ be nonnegative functions in  $L^{\infty}(\partial \Omega),$ given   $s\geq 1 ,$   $0<\gamma\leq1$ and $ \eta>0.$ Then there exists a weak solution to \eqref{1}.
\end{thm}
\subsection{Approximation Scheme}
In this section, to demonstrate Theorem \ref{48}, we will focus on the following approximating problem
\begin{equation}\label{9}
    \left\{
  \begin{array}{ll}
   -\mbox{div }(\aleph(y, \nabla u_{n}))=\frac{f_{n}}{(\vert u_{n}\vert+\frac{1}{n})^{\gamma}}& \hbox{in $\Omega,$} \\
  \frac{\partial u_{n}}{\partial \sigma}+\lambda\vert u_{n}\vert^{p-2} u_{n}+\vert u_{n}\vert^{s-1} u_{n}=\frac{g_{n}}{(\vert u_{n}\vert+\frac{1}{n})^{\eta}} & \hbox{on $\partial\Omega,$} \\
  \end{array}
\right.
\end{equation}
with $ f_{n}=T_{n}(f)$ and $ g_{n}=T_{n}(g).$ We define a weak solution $ u_{n}\in W^{1,p}(\Omega)$ of \eqref{9} satisfying
\begin{align*}
&\int_{\Omega}\aleph(y,\nabla u_{n})\nabla \varphi dy+\int_{\partial \Omega} \lambda \vert u_{n}\vert^{p-2}u_{n}\varphi d\mathcal{H}^{m-1}+\int_{\partial \Omega}\vert u_{n}\vert^{s-1}u_{n}\varphi d\mathcal{H}^{m-1}\nonumber \\
&=\int_{\Omega}\frac{f_{n}}{(\vert u_{n}\vert+\frac{1}{n})^{\gamma}}\varphi dy+\int_{\partial \Omega}\frac{g_{n}}{(\vert u_{n}\vert+\frac{1}{n})^{\eta}}\varphi d\mathcal{H}^{m-1},\, \forall \varphi \in W^{1,p}(\Omega)\cap L^{\infty}(\partial \Omega).
\end{align*}
\par Let's note that, for all that follows, the $C_i$, where $i = 1,..,$ represent positive constants. Throughout the article, a "test function" will be abbreviated as \textit{t. f.}.
\begin{lem}
Assume that  $0\leq f\in L^{\frac{m}{p}}(\Omega),$ let  $ \lambda$ and $ g$ be nonnegative functions in  $L^{\infty}(\partial \Omega),$ given $ s\geq 1,$ $0<\gamma\leq1$ and $ \eta>0.$ Then problem \eqref{9} admits a nonnegative solution $ u_{n}$ in $ W^{1,p}(\Omega).$
\end{lem}
\begin{dem}
Let $ n\in \mathbb{N}$ be fixed and suppose  $\vartheta\in L^{p}(\partial\Omega).$ Now, let's  consider the following problem
\begin{equation}\label{2}
    \left\{
  \begin{array}{ll}
   -\textnormal{div }(\aleph(y, \nabla w))=\frac{f_{n}}{(\vert w\vert+\frac{1}{n})^{\gamma}}& \hbox{in $\Omega,$} \\
  \frac{\partial w}{\partial \sigma}+\lambda\vert w\vert^{p-2}w+\vert w\vert^{s-1}w=\frac{g_{n}}{(\vert \vartheta\vert+\frac{1}{n})^{\eta}} & \hbox{on $\partial\Omega.$} \\
  \end{array}
\right.
\end{equation}
The existence of a solution $ w\in W^{1,p}(\Omega)$ to \eqref{2} is assured by \textnormal{\cite{7030}}.
\par Taking $w^{-}$ as a \textit{t. f.} in \eqref{2}  we get
\begin{equation*}
-\int_{\Omega} \aleph(y, \nabla w)\nabla w dy-\int_{\partial \Omega}\lambda \vert w^{-}\vert^{p}d\mathcal{H}^{m-1}-\int_{\partial \Omega}\vert w^{-}\vert^{s+1}d\mathcal{H}^{m-1}=\int_{\Omega}\frac{f_{n}w^{-}}{(\vert w\vert+\frac{1}{n})^{\gamma}}dy+\int_{\partial \Omega}\frac{g_{n}w^{-}}{(\vert \vartheta\vert+\frac{1}{n})^{\eta}}d\mathcal{H}^{m-1},
\end{equation*}
which gives  $w^{-}=0 ,$ hence $ w\geq 0$ $ \mathcal{H}^{m-1}$ a. e. (almost everywhere) on $\partial \Omega$ and a. e.  in $ \Omega.$
\par Now, we will demonstrate the existence of a fixed of  the map
$S:L^{p}(\partial \Omega)\rightarrow L^{p}(\partial \Omega) $ with $S(\vartheta)=w|_{\partial \Omega}.$ Since the datum  $\frac{g_{n}}{(\vert \vartheta\vert+\frac{1}{n})^{\eta}}$ is bounded, there exists $ C_{2}>0$ independent of $ \vartheta$ and $ w$ but possibly depending on $ n,$ such that $\Vert w \Vert_{L^{\infty}(\partial \Omega)}\leq C_{2}.$ We choose  $w$ as a \textit{t. f.} in \eqref{2} and by \eqref{b}, it yields
\begin{align*}
\alpha\int_{\Omega}\vert \nabla w\vert^{p} dy+\int_{\partial \Omega}\lambda \vert w\vert^{p}d\mathcal{H}^{m-1}+\int_{\partial \Omega}\vert w \vert^{s+1}d\mathcal{H}^{m-1}&\leq \int_{\Omega}\frac{f_{n}w}{(\vert w\vert+\frac{1}{n})^{\gamma}} dy+\int_{\partial\Omega}\frac{g_{n}w}{(\vert \vartheta \vert+\frac{1}{n})^{\eta}} d\mathcal{H}^{m-1}\\
&\leq n^{\gamma+1}\int_{\Omega} \vert w \vert dy+n^{\eta+1}\int_{\partial\Omega}\vert w \vert d\mathcal{H}^{m-1}.
\end{align*}
Applying Young's inequality  we obtain
\begin{align}\label{125}
&\vert\vert w\vert \vert^{p}_{\frac{\lambda}{\alpha},p}\leq n^{\gamma+1}\varepsilon^{(p-1)p}_{1}\int_{\Omega}\vert w\vert^{p}dy+n^{\eta+1}\varepsilon^{(p-1)p}_{2}\int_{\partial \Omega}\vert w\vert^{p}d\mathcal{H}^{m-1}+n^{\gamma+1}\varepsilon^{-p}_{1}\vert \Omega\vert+n^{\eta+1}\varepsilon^{-p}_{2}\mathcal{H}^{m-1}(\partial \Omega),
\end{align}
where $\varepsilon_{1}$ and $\varepsilon_{2}$ are any positive constants. From \eqref{46} and the fact that $\Vert.\Vert_{\frac{\lambda}{\alpha},p}$ and $ \Vert.\Vert_{W^{1,p}(\Omega)}$  are equivalent norms, we conclude that
\begin{equation*}
    \int_{\partial \Omega}w^{p}d\mathcal{H}^{m-1}\leq \bigg(\int_{\partial \Omega}\vert w\vert^{\frac{(m-1)p}{m-p}}d\mathcal{H}^{m-1}\bigg)^{\frac{m-p}{m-1}}\bigg(\mathcal{H}^{m-1}(\partial \Omega)\bigg)^{\frac{p-1}{m-1}}\leq C_{3}\Vert w\Vert^{p}_{\frac{\lambda}{\alpha},p}.
\end{equation*}
Using previous inequality, then \eqref{125} becomes
\begin{align*}
&\vert\vert w\vert \vert^{p}_{\frac{\lambda}{\alpha},p}\leq n^{\gamma+1}\varepsilon^{(p-1)p}_{1}\Vert w\Vert_{\frac{\lambda}{\alpha},p}^{p} +n^{\eta+1}\varepsilon^{(p-1)p}_{2}C_{3}\Vert w\Vert_{\frac{\lambda}{\alpha},p}^{p}+n^{\gamma+1}\varepsilon^{-p}_{1}\vert \Omega\vert+n^{\eta+1}\varepsilon^{-p}_{2}\mathcal{H}^{m-1}(\partial \Omega),
\end{align*}
choosing $ \varepsilon_{1}$ satisfying   $n^{\gamma+1}\varepsilon^{(p-1)p}_{1}<\frac{1}{4}$ and $ \varepsilon_{2}$ such that $n^{\eta+1}\varepsilon^{(p-1)p}_{2}C_{3}<\frac{1}{4},$ we deduce that
\begin{align*}
& \Vert w\Vert_{\frac{\lambda}{\alpha},p}\leq  C_{4}.
\end{align*}
We know that  $\Vert.\Vert_{\frac{\lambda}{\alpha},p}$ and $ \Vert.\Vert_{W^{1,p}(\Omega)}$ are equivalent norms, then
\begin{align}\label{124}
& \Vert w\Vert_{W^{1,p}(\Omega)}\leq  C_{5}.
\end{align}
By \eqref{46} and \eqref{124} we have that
\begin{align*}
& \Vert w \Vert_{L^{p}(\partial \Omega)} \leq  C_{1}\vert\vert w\vert\vert_{W^{1,p}(\Omega)}(\mathcal{H}^{m-1}(\Omega))^{\frac{p-1}{p(m-1)}}\leq C_{6},
\end{align*}
with $C_{6}$ independent of $ w.$ Therefore,  the ball $B $ in $L^{p}(\partial\Omega)$ with radius  $ C_{6}$ is invariant for $S.$ Moreover by the compactness of the trace embedding  $S(B) $ is relatively compact in $ L^{p}(\partial \Omega).$ We shall now prove that, on $B $, the function $ S$ is continuous. Assume that the sequence of functions $\vartheta_{k}$ in the ball $B$ converges to $\vartheta$ in $L^{p}(\partial\Omega)$ as $k$ approaches $+\infty,$ and we consider
$ S(\vartheta_{k})=w_{k}|_{\partial \Omega} $ such that
\begin{equation}\label{8}
    \left\{
  \begin{array}{ll}
   -\textnormal{div }(\aleph(y, \nabla w_{k}))=\frac{f_{n}}{(\vert w_{k}\vert+\frac{1}{n})^{\gamma}}& \hbox{in $\Omega,$} \\
  \frac{\partial w_{k}}{\partial \sigma}+\lambda\vert w_{k}\vert^{p-2}w_{k}+\vert w_{k}\vert^{s-1}w_{k}=\frac{g_{n}}{(\vert \vartheta_{k}\vert+\frac{1}{n})^{\eta}} & \hbox{on $\partial\Omega.$} \\
  \end{array}
\right.
\end{equation}
We take  $ w_{k}$ as a \textit{t. f.} in \eqref{8} and using \eqref{b}, we get
\begin{equation*}
\alpha\int_{\Omega}\vert \nabla w_{k}\vert^{p} dy +\int_{\partial \Omega}\lambda \vert w_{k} \vert^{p} d\mathcal{H}^{m-1}+\int_{\partial \Omega}\vert w_{k} \vert^{s+1} d\mathcal{H}^{m-1}\leq \int_{\Omega}\frac{f_{n}w_{k}}{(\vert w_{k}\vert+\frac{1}{n})^{\gamma}}dy+\int_{\partial \Omega}\frac{g_{n}w_{k}}{(\vert \vartheta_{k} \vert+\frac{1}{n})^{\eta}}d\mathcal{H}^{m-1}.
\end{equation*}
From Young's inequality and \eqref{46}, we deduce that
\begin{align}
&\Vert w_{k} \Vert_{W^{1,p}(\Omega)}\leq C_{7},\label{56}\\
&\Vert w_{k} \Vert_{L^{p}(\partial \Omega)}\leq C_{8},\\
&\Vert w_{k} \Vert_{L^{s+1}(\partial \Omega)}\leq C_{9},\\ 
\end{align}
where $C_{7},$ $ C_{8}$ and $ C_{9}$ are independent of $k.$
\par Let $ \tau>0,$ we take $ T_{\tau}(w_{k})-T_{\tau}(w)$ as \textit{t. f.} in \eqref{8}, we get
\begin{align*}
&\int_{\Omega}\aleph(y,\nabla w_{k})\nabla(T_{\tau}(w_{k})-T_{\tau}(w))dy+\int_{\partial \Omega}\vert w_{k}\vert^{p-1}(T_{\tau}(w_{k})-T_{\tau}(w))d\mathcal{H}^{m-1}+\int_{\partial \Omega}\vert w_{k}\vert^{s-1}(T_{\tau}(w_{k})-T_{\tau}(w))d\mathcal{H}^{m-1}\nonumber\\
&=\int_{\Omega}\frac{f_{n}}{(\vert w_{k}\vert+\frac{1}{n})^{\gamma}}(T_{\tau}(w_{k})-T_{\tau}(w))dy+\int_{\partial \Omega}\frac{g_{n}}{(\vert \vartheta_{k} \vert+\frac{1}{n})^{\eta}}(T_{\tau}(w_{k})-T_{\tau}(w))d\mathcal{H}^{m-1}.
\end{align*}
By \eqref{56}, \eqref{46} and Lebesgue's Theorem, we obtain
\begin{equation}\label{131}
\lim_{k\rightarrow +\infty}\int_{\partial \Omega}\vert w_{k}\vert^{p-1}(T_{\tau}(w_{k})-T_{\tau}(w))d\mathcal{H}^{m-1}=0,
\end{equation}
\begin{equation}\label{132}
\lim_{k\rightarrow +\infty}\int_{\partial \Omega}\vert w_{k}\vert^{s-1}(T_{\tau}(w_{k})-T_{\tau}(w))d\mathcal{H}^{m-1}=0,
\end{equation}
and
\begin{equation}\label{133}
 \lim_{k\rightarrow +\infty}\int_{\partial \Omega}\frac{g_{n}}{(\vert \vartheta_{k} \vert+\frac{1}{n})^{\eta}}(T_{\tau}(w_{k})-T_{\tau}(w))d\mathcal{H}^{m-1}=0.
\end{equation}
Using \eqref{56} and Lebesgue's Theorem, we have
\begin{equation}\label{134}
 \lim_{k\rightarrow +\infty}\int_{\Omega}\frac{f_{n}}{(\vert w_{k}\vert+\frac{1}{n})^{\gamma}}(T_{\tau}(w_{k})-T_{\tau}(w))dy=0.
\end{equation}
We see that
\begin{equation}\label{59}
 \int_{\Omega}\aleph(y,\nabla w_{k}).\nabla(T_{\tau}(w_{k})-T_{\tau}(w))dy=\int_{\Omega}\aleph(y,\nabla T_{\tau}(w_{k})).\nabla(T_{\tau}(w_{k})-T_{\tau}(w))dy-\int_{\{w_{k}>\tau\}}\aleph(y,\nabla w_{k}).\nabla T_{\tau}(w)dy.
\end{equation}
We have that  $\aleph(y,\nabla w_{k})$ is bounded in $L^{p^{\prime}}(\Omega)$ and that $\chi_{\{w_{k}>\tau\}}\nabla T_{\tau}(w)\rightarrow 0$ strongly in $L^{p}(\Omega).$ Then
\begin{align}\label{60}
&\lim_{k\rightarrow +\infty}\int_{\{w_{k}>\tau\}}\aleph(y,\nabla w_{k}).\nabla T_{\tau}(w)dy=0.
\end{align}
We can formulate that
\begin{align}\label{110}
\int_{\Omega}\aleph(y,\nabla T_{\tau}(w_{k})).\nabla(T_{\tau}(w_{k})-T_{\tau}(w))dy
&=\int_{\Omega}(\aleph(y,\nabla T_{\tau}(w_{k}))-\aleph(y,\nabla T_{\tau}(w))).\nabla(T_{\tau}(w_{k})-T_{\tau}(w))dy\nonumber\\
&+\int_{\Omega}\aleph(y,\nabla T_{\tau}(w)).\nabla(T_{\tau}(w_{k})-T_{\tau}(w))dy.
\end{align}
Keeping  in mind that  $ \aleph(y,\nabla T_{\tau}(w))\in L^{p'}(\Omega)$ and that $ T_{\tau}(w_{k}) $ converges to $ T_{\tau}(w) $ weakly in $W^{1,p}(\Omega),$ we find
\begin{align}\label{130}
&\lim_{k\rightarrow +\infty}\int_{\Omega}\aleph(y,\nabla T_{\tau}(w)).\nabla(T_{\tau}(w_{k})-T_{\tau}(w))dy=0.
\end{align}
By \eqref{59}, \eqref{60}, \eqref{110} and  \eqref{130}, we see
\begin{equation*}
\lim_{k\rightarrow +\infty}\int_{\Omega}\aleph(y,\nabla w_{k}).\nabla(T_{\tau}(w_{k})-T_{\tau}(w))dy =\lim_{k\rightarrow +\infty}\int_{\Omega}(\aleph(y,\nabla T_{\tau}(w_{k}))-\aleph(y,\nabla T_{\tau}(w))).\nabla(T_{\tau}(w_{k})-T_{\tau}(w))dy.
\end{equation*}
Utilizing \eqref{131}, \eqref{132}, \eqref{133}, and \eqref{134}, we obtain
\begin{equation*}
\lim_{k\rightarrow +\infty}\int_{\Omega}(\aleph(y,\nabla T_{\tau}(w_{k}))-\aleph(y,\nabla T_{\tau}(w))).\nabla(T_{\tau}(w_{k})-T_{\tau}(w))dy=0,
\end{equation*}
which permits us to use \cite[Lemma 5]{6014}, enabling us to conclude
\begin{equation*}
    T_{\tau}(w_{k})\rightarrow T_{\tau}(w) \hspace*{0.2cm} strongly \hspace*{0.1cm} in \hspace*{0.1cm} W^{1,p}(\Omega).
\end{equation*}
This is sufficient to conclude $ w_{k}=S(\vartheta_{k}) \rightarrow w=S(\vartheta)$ in $ L^{p}(\partial \Omega).$ Then $ S$ is continuous on $ B.$ Due to the Schauder fixed point Theorem, a solution $ u_{n}\in W^{1,p}(\Omega)\cap L^{\infty}(\partial \Omega)$ exists for \eqref{9}.
\begin{flushright}
    $\square$
\end{flushright}
\end{dem}
\begin{lem}
Suppose $ f\in L^{\frac{m}{p}}(\Omega)$ is a nonnegative function, let $ \lambda$ and $ g$ be nonnegative  functions in  $L^{\infty}(\partial \Omega),$ given $s\geq 1 ,$   $0<\gamma\leq1$ and $ \eta>0,$ let $ u_{n}$ be a solution to \eqref{9}. Hence there exists $ \beta_{p}>0$ independent of $n$  such that
\begin{equation*}
\forall n\in\mathbb{N}, \hspace*{0.5cm} u_{n}\geq \beta_{p}>0 \textnormal{ for } \mathcal{H}^{m-1} \textnormal{ a. e.  on } \partial \Omega.
\end{equation*}
\end{lem}
\begin{dem}
Choosing  $(u_{n}-u_{n+1})^{+}$ as a \textit{t. f.} in the difference between the problems solved respectively by $ u_{n}$ and $ u_{n+1},$ using the facts $f_{n}\leq f_{n+1}$ and $g_{n}\leq g_{n+1},$ we find that
\begin{align*}
&\int_{\Omega}(\aleph(y,\nabla u_{n})-\aleph(y,\nabla u_{n+1})) \nabla(u_{n}-u_{n+1})^{+}dy\\
&+\int_{\partial\Omega}\lambda (\vert u_{n}\vert^{p-1}-\vert u_{n+1}\vert^{p-1})(u_{n}-u_{n+1})^{+}d\mathcal{H}^{m-1}+\int_{\partial\Omega}(\vert u_{n}\vert^{s}-\vert u_{n+1}\vert^{s})(u_{n}-u_{n+1})^{+}d\mathcal{H}^{m-1}\\
&=\int_{\Omega}(\frac{f_{n}}{(u_{n}+\frac{1}{n})^{\gamma}}-\frac{f_{n+1}}{(u_{n+1}+\frac{1}{n+1})^{\gamma}})(u_{n}-u_{n+1})^{+}dy+\int_{\partial\Omega}(\frac{g_{n}}{(u_{n}+\frac{1}{n})^{\gamma}}-\frac{g_{n+1}}{(u_{n+1}+\frac{1}{n+1})^{\gamma}})(u_{n}-u_{n+1})^{+}d\mathcal{H}^{m-1}\\
&\leq\int_{\Omega}\frac{f_{n+1}((u_{n+1}+\frac{1}{n+1})^{\gamma}-(u_{n}+\frac{1}{n})^{\gamma})}{(u_{n}+\frac{1}{n})^{\gamma}(u_{n+1}+\frac{1}{n+1})^{\gamma}}(u_{n}-u_{n+1})^{+}dy+\int_{\partial\Omega}\frac{g_{n+1}((u_{n+1}+\frac{1}{n+1})^{\eta}-(u_{n}+\frac{1}{n})^{\eta})}{(u_{n}+\frac{1}{n})^{\eta}(u_{n+1}+\frac{1}{n+1})^{\eta}}(u_{n}-u_{n+1})^{+}d\mathcal{H}^{m-1}\leq0.
\end{align*}
By \eqref{d}, we have
\begin{align*}
& 0\leq\int_{\Omega}(\aleph(y,\nabla u_{n})-\aleph(y,\nabla u_{n+1})) \nabla(u_{n}-u_{n+1})^{+}dy.
\end{align*}
Moreover,
\begin{equation*}
0\leq\int_{\partial\Omega}\lambda(\vert u_{n}\vert^{p-1}-\vert u_{n+1}\vert^{p-1})(u_{n}-u_{n+1})^{+}d\mathcal{H}^{m-1}
\mbox{ and }
0\leq\int_{\partial\Omega}(\vert u_{n}\vert^{s}-\vert u_{n+1}\vert^{s})(u_{n}-u_{n+1})^{+}d\mathcal{H}^{m-1}.
\end{equation*}
Therefore, $(u_{n}-u_{n+1})^{+}=0,$ leading to  $  u_{n}\leq u_{n+1}$ $ \mathcal{H}^{m-1}$ a. e. on $\partial \Omega$ and a. e. in $ \Omega.$
We observe the existence of a nonnegative solution $\vartheta\in C^{1}(\overline{\Omega})$ to
\begin{equation}\label{33}
    \left\{
  \begin{array}{ll}
-\textnormal{div }(\aleph(y,\nabla \vartheta))=\frac{f_{1}}{(\vartheta+1)^{\gamma}}& \hbox{ in } \Omega, \\
\frac{\partial \vartheta}{\partial \sigma}+\lambda \vartheta^{p-1}+\vartheta^{s}=0  & \hbox{ on } \partial\Omega. \\
  \end{array}
\right.
\end{equation}
By \cite[Theorem 2]{6025} we deduce that  $ \vartheta>0$ in $ \overline{\Omega}.$ Subtracting \eqref{9} from \eqref{33} we get
\begin{equation}\label{34}
    \left\{
  \begin{array}{ll}
-\textnormal{div }\big(\aleph(y,\nabla u_{n})-\aleph(y,\nabla \vartheta)\big)=\frac{f_{n}}{(\vert u_{n}\vert+\frac{1}{n})^{\gamma}}-\frac{f_{1}}{(\vartheta+1)^{\gamma}}& \hbox{in $\Omega,$} \\
\frac{\partial(u_{n}- \vartheta)}{\partial \sigma}+\lambda(u_{n}^{p-1}-\vartheta^{p-1})+u_{n}^{s}-\vartheta^{s}=\frac{g_{n}}{(\vert u_{n}\vert+\frac{1}{n})^{\eta}}  & \hbox{ on $\partial\Omega.$} \\
  \end{array}
\right.
\end{equation}
We take  $(\vartheta-u_{n})^{+}$ as a \textit{t. f.} in problem \eqref{34}, we can prove that
\begin{equation*}
\forall n\in\mathbb{N},\hspace*{0.4cm} u_{n}\geq \vartheta \textnormal{ for } \mathcal{H}^{m-1} \textnormal{a. e. on } \partial \Omega.
\end{equation*}
We know that $ \vartheta$ is continuous with $ \vartheta>0$ on $\partial \Omega,$ then
\begin{equation*}
\forall n\in\mathbb{N}, \hspace*{0.4cm}u_{n}\geq \vartheta>\min_{\partial \Omega}\vartheta=\beta_{p} \hspace*{0.1cm}\textnormal{for}\hspace*{0.1cm} \mathcal{H}^{m-1} \textnormal{ a. e. on } \partial \Omega.
\end{equation*}
\begin{flushright}
    $\square$
\end{flushright}
\end{dem}
\begin{lem}\label{85}
Supposing that $ f\in L^{\frac{m}{p}}(\Omega)$ is a nonnegative function, let  $ \lambda$ and $ g$ are nonnegative functions in  $L^{\infty}(\partial \Omega),$ given $s\geq 1,$   $0<\gamma\leq1$ and $ \eta>0.$ Assume that $ u_{n}$ is a solution to problem \eqref{9}. Then
\begin{equation}\label{40}
    \Vert u_{n}\Vert_{W^{1,p}(\Omega)}\leq C_{10},
\end{equation}
and
\begin{align}\label{43}
\int_{\partial \Omega}\vert u_{n}\vert^{s+1}d\mathcal{H}^{m-1}\leq C_{11},
\end{align}
with the constants $ C_{10}$ and $C_{11}$ are independent of $ n.$  Moreover there exists a subsequence $ u_{n}$ and a function $ u_{p}$  such that
\begin{align}
& u_{n}\rightharpoonup u_{p}\hspace*{0.2cm}\textnormal{weakly}\hspace*{0.2cm}\textnormal{in}\hspace*{0.2cm}W^{1,p}(\Omega),\nonumber\\
& u_{n}\rightarrow u_{p}\hspace*{0.2cm}\textnormal{strongly}\hspace*{0.2cm}\textnormal{in}\hspace*{0.2cm} L^{q}(\Omega)\hspace*{0.2cm} for\hspace*{0.1cm}any \hspace*{0.1cm}  1\leq q<\frac{pm}{m-p},  \nonumber \\ 
& u_{n}\rightarrow u_{p}\hspace*{0.2cm}\textnormal{strongly}\hspace*{0.2cm}\textnormal{in}\hspace*{0.2cm} L^{t}(\partial\Omega)\hspace*{0.2cm} for\hspace*{0.1cm}any \hspace*{0.1cm}  1<t<\frac{p(m-1)}{m-p},\label{127}\\
&\vert u_{n}\vert^{s-1}u_{n}\rightharpoonup \vert u_{p}\vert^{s-1}u_{p}\hspace*{0.2cm}\textnormal{weakly}\hspace*{0.2cm}\textnormal{in}\hspace*{0.2cm}L^{\frac{s+1}{s}}(\partial \Omega).\label{38}
\end{align}
\end{lem}
\begin{dem}
We take  $ u_{n}$ as a \textit{t. f.} in \eqref{9} and applying \eqref{b}, we obtain
\begin{equation}\label{32}
\alpha\int_{\Omega}\vert\nabla u_{n} \vert^{p}dy+ \int_{\partial \Omega}\lambda \vert u_{n} \vert^{p}d\mathcal{H}^{m-1}+\int_{\partial \Omega}\vert u_{n}\vert^{s+1}d\mathcal{H}^{m-1}\leq\int_{\Omega}\frac{f_{n}u_{n}}{(u_{n}+\frac{1}{n})^{\gamma}}dy+\int_{\partial\Omega}\frac{g_{n}u_{n}}{(u_{n}+\frac{1}{n})^{\eta}}d\mathcal{H}^{m-1}.
\end{equation}
For the first integral in \eqref{32} on the right-hand side in the case where $ 0<\gamma<1,$ employing H\"{o}lder's and Sobolev's inequalities we have
\begin{equation}\label{62}
\int_{\Omega}\frac{f_{n}u_{n}}{(u_{n}+\frac{1}{n})^{\gamma}}dy\leq \int_{\Omega} f_{n}u^{1-\gamma}_{n} dy\leq \Vert f_{n}\Vert_{L^{\frac{m}{p}}(\Omega)}\vert\Omega\vert^{\frac{p-1+\gamma}{p}(\frac{m-p}{m})}\Vert u_{n}\Vert^{1-\gamma}_{L^{\frac{mp}{m-p}}(\Omega)}\leq \Vert f \Vert_{L^{\frac{m}{p}}(\Omega)}\vert\Omega\vert^{\frac{p-1+\gamma}{p}(\frac{m-p}{m})}C_{m,p}^{1-\gamma}\Vert u_{n}\Vert^{1-\gamma}_{W^{1,p}(\Omega)},
\end{equation}
with $ C_{m,p}$ is the best constant in Sobolev's inequality. If $ \gamma=1,$ we have
\begin{align}\label{98}
&\int_{\Omega}\frac{f_{n}u_{n}}{u_{n}+\frac{1}{n}}dy\leq \int_{\Omega}fdy.
\end{align}
For the second  integral on the right hand side of \eqref{32} in the case where $ \eta <1,$ utilizing H\"{o}lder's inequality and \eqref{46} we see
\begin{align}\label{63}
&\int_{\partial\Omega}\frac{g_{n}u_{n}}{(u_{n}+\frac{1}{n})^{\eta}}d\mathcal{H}^{m-1}\leq\int_{\Omega}g_{n}u^{1-\eta}_{n}dy\leq
\Vert g_{n}\Vert_{L^{\infty}(\partial \Omega)}\Vert u_{n}\Vert^{1-\eta}_{W^{1,p}(\Omega)}\big(\mathcal{H}^{m-1}(\partial \Omega)\big)^{\frac{p(m-1)-(1-\eta)(m-p)}{p(m-1)}}.
\end{align}
If $ \eta\geq 1$ we arrive at
\begin{align}\label{97}
&\int_{\partial\Omega}\frac{g_{n}u_{n}}{(u_{n}+\frac{1}{n})^{\eta}}d\mathcal{H}^{m-1}\leq\frac{\Vert g_{n}\Vert_{L^{\infty}(\partial \Omega)}}{\beta_{p}^{\eta-1}}.
\end{align}
By \eqref{32}, \eqref{62} (respectively \eqref{98})  and \eqref{63} (respectively \eqref{97}) we can write
\begin{align}\label{61}
\alpha\int_{\Omega}\vert\nabla u_{n} \vert^{p}dy+ \int_{\partial \Omega}\lambda \vert u_{n} \vert^{p}d\mathcal{H}^{m-1}+\int_{\partial \Omega}\vert u_{n} \vert^{s+1}d\mathcal{H}^{m-1}&\leq \Vert f \Vert_{L^{\frac{m}{p}}(\Omega)}\vert\Omega\vert^{\frac{p-1+\gamma}{p}(\frac{m-p}{m})}C_{m,p}^{1-\gamma}\Vert u_{n}\Vert^{1-\gamma}_{W^{1,p}(\Omega)}\nonumber\\
&+\Vert g\Vert_{L^{\infty}(\partial \Omega)}\Vert u_{n}\Vert^{1-\eta}_{W^{1,p}(\Omega)}\big(\mathcal{H}^{m-1}(\partial \Omega)\big)^{\frac{p(m-1)-(1-\eta)(m-p)}{p(m-1)}}.
\end{align}
Having in mind that  $ \Vert.\Vert_{\frac{\lambda}{\alpha},p}$ and $\Vert .\Vert_{W^{1,p}(\Omega)} $ are equivalent norms and  using Young's inequality in the right hand side of \eqref{61} we have
\begin{align*}
\Vert u_{n}\Vert_{W^{1,p}(\Omega)} \leq C_{10}.
\end{align*}
Hence there exists a subsequence $ u_{n}$ and a function $ u_{p}$  such that
\begin{align*}
& u_{n}\rightharpoonup u_{p}\hspace*{0.2cm}\textnormal{weakly}\hspace*{0.2cm}\textnormal{in}\hspace*{0.2cm}W^{1,p}(\Omega),\\
& u_{n}\rightarrow u_{p}\hspace*{0.2cm}\textnormal{strongly}\hspace*{0.2cm}\textnormal{in}\hspace*{0.2cm} L^{q}(\Omega)\hspace*{0.2cm} for\hspace*{0.1cm}any \hspace*{0.1cm}  1\leq q<\frac{pm}{m-p},\\
& u_{n}\rightarrow u_{p}\hspace*{0.2cm}\textnormal{strongly}\hspace*{0.2cm}\textnormal{in}\hspace*{0.2cm} L^{t}(\partial\Omega)\hspace*{0.2cm} for\hspace*{0.1cm}any \hspace*{0.1cm}  1<t<\frac{p(m-1)}{m-p}.
\end{align*}
Moreover we observe that
\begin{align*}
&\int_{\partial \Omega}(\vert u_{n}\vert^{s})^{\frac{s+1}{s}}d\mathcal{H}^{m-1}=\int_{\partial \Omega}\vert u_{n}\vert^{s+1}d\mathcal{H}^{m-1}\leq C_{11}.
\end{align*}
Then
\begin{align*}
&\vert u_{n}\vert^{s-1}u_{n}\rightharpoonup \vert u_{p}\vert^{s-1}u_{p}\hspace*{0.2cm}\textnormal{weakly}\hspace*{0.2cm}\textnormal{in}\hspace*{0.2cm}L^{\frac{s+1}{s}}(\partial \Omega).
\end{align*}
\begin{flushright}
    $\square$
\end{flushright}
\end{dem}
\begin{lem}\label{99}
Let  $ f\in L^{\frac{m}{p}}(\Omega)$ is nonnegative function, let  $ \lambda$ and $ g$ are nonnegative functions in  $L^{\infty}(\partial \Omega),$ supposing that   $s\geq 1 ,$ $0<\gamma\leq1$ and $ \eta>0,$ let $ u_{p}$ be a function found in Lemma \ref{85}. Then
\begin{equation*}
    \frac{f}{u^{\gamma}_{p}}\in L^{1}_{loc}(\Omega).
\end{equation*}
\end{lem}
\begin{dem}
In \eqref{1}, taking $\varphi \in W^{1,p}(\Omega)\cap L^{\infty}(\Omega)$ as \textit{t. f.},  we obtain
\begin{align}\label{86}
&\int_{\Omega}\aleph(y,\nabla u_{n})\nabla \varphi dy+\int_{\partial \Omega} \lambda \vert u_{n}\vert^{p-2}u_{n}\varphi d\mathcal{H}^{m-1}+\int_{\partial \Omega}\vert u_{n}\vert^{s-1}u_{n}\varphi d\mathcal{H}^{m-1}\nonumber \\
&=\int_{\Omega}\frac{f_{n}}{(u_{n}+\frac{1}{n})^{\gamma}}\varphi dy+\int_{\partial \Omega}\frac{g_{n}}{(u_{n}+\frac{1}{n})^{\eta}}\varphi d\mathcal{H}^{m-1}.
\end{align}
Dropping nonnegative term  on the right  hand side of \eqref{86}, by \eqref{c} and  Young's inequality we find  that
\begin{align*}
\int_{\Omega}\frac{f_{n}}{(u_{n}+\frac{1}{n})^{\gamma}}\varphi dy &\leq\int_{\Omega}\aleph(y,\nabla u_{n})\nabla \varphi dy+\int_{\partial \Omega} \lambda \vert u_{n}\vert^{p-1}\varphi d\mathcal{H}^{m-1}+\int_{\partial \Omega}\vert u_{n}\vert^{s}\varphi d\mathcal{H}^{m-1}\\
&\leq\int_{\Omega}\vert l\vert \vert \nabla\varphi\vert  dy+\frac{p-1}{p}\int_{\Omega}\vert\nabla u_{n}\vert^{p}dy+\frac{1}{p}\int_{\Omega}\vert\nabla \varphi\vert^{p}dy+\frac{p-1}{p}\int_{\partial\Omega}\lambda \vert u_{n}\vert^{p} d\mathcal{H}^{m-1}\\
&+\frac{1}{p}\int_{\partial \Omega}\lambda\vert \varphi \vert^{p} d\mathcal{H}^{m-1}+\frac{s}{s+1}\int_{\partial \Omega}\vert u_{n}\vert^{s+1}d\mathcal{H}^{m-1}+\frac{1}{s+1}\int_{\partial \Omega}\vert \varphi\vert^{s+1}d\mathcal{H}^{m-1}\\
&\leq C_{12},
\end{align*}
where $C_{12}$ is independent of $n.$ Fatou's Lemma then implies
\begin{equation}\label{41}
    \int_{\Omega}\frac{f}{u_{p}^{\gamma}}\varphi dy\leq C_{12}.
\end{equation}
This implies that
\begin{equation*}
    \frac{f}{u^{\gamma}_{p}}\in L^{1}_{loc}(\Omega).
\end{equation*}
\begin{flushright}
    $\square$
\end{flushright}
\end{dem}
\begin{lem}\label{96}
Let  $ f\in L^{\frac{m}{p}}(\Omega)$ is nonnegative function, let  $ \lambda$ and $ g$ are nonnegative functions in  $L^{\infty}(\partial \Omega),$ supposing that    $s\geq 1 ,$   $0<\gamma\leq1$ and $ \eta>0,$ let $ u_{n}$ be a solution to \eqref{9}. Then $T_{\tau}(u_{n})$ converge strongly to $T_{\tau}(u_{p})$ in $ W^{1,p}(\Omega).$
\end{lem}
\begin{dem}
Taking $(T_{\tau}(u_{n})-T_{\tau}(u_{p}))$ as a \textit{t. f.} in \eqref{9},  we  find that
\begin{align}\label{77}
&\int_{\Omega}\aleph(y,\nabla u_{n})\nabla(T_{\tau}(u_{n})-T_{\tau}(u_{p}))dy+\int_{\partial \Omega}\lambda\vert u_{n}\vert^{p-2}u_{n}(T_{\tau}(u_{n})-T_{\tau}(u_{p}))d\mathcal{H}^{m-1}\nonumber\\
&+\int_{\partial \Omega}\vert u_{n}\vert^{s-1}\, u_{n}(T_{\tau}(u_{n})-T_{\tau}(u_{p}))d\mathcal{H}^{m-1}=\int_{\Omega}\frac{f_{n}}{(u_{n}+\frac{1}{n})^{\gamma}}(T_{\tau}(u_{n})-T_{\tau}(u_{p}))dy\nonumber\\
&+\int_{\partial \Omega}\frac{g_{n}}{(u_{n}+\frac{1}{n})^{\eta}}(T_{\tau}(u_{n})-T_{\tau}(u_{p}))d\mathcal{H}^{m-1}.
\end{align}
For the first term  on the left  hand side of \eqref{77}, we can write
\begin{equation}\label{78}
 \int_{\Omega}\aleph(y,\nabla u_{n}).\nabla(T_{\tau}(u_{n})-T_{\tau}(u_{p}))dy=\int_{\Omega}\aleph(y,\nabla T_{\tau}(u_{n})).\nabla(T_{\tau}(u_{n})-T_{\tau}(u_{p}))dy-\int_{\{u_{n}>\tau\}}\aleph(y,\nabla u_{n}).\nabla T_{\tau}(u_{p})dy.
\end{equation}
Since $\aleph(y,\nabla u_{n})$ is bounded in $L^{p^{\prime}}(\Omega)$ and $\chi_{\{u_{n}>\tau\}}\nabla T_{\tau}(u_{p})\rightarrow 0$ strongly in $L^{p}(\Omega),$ then
\begin{align*}
&\lim_{n\rightarrow +\infty}\int_{\{u_{n}>\tau\}}\aleph(y,\nabla u_{n}).\nabla T_{\tau}(u_{p})dy=0.
\end{align*}
For the first term on the right hand side of \eqref{78} we can estimate
\begin{align}\label{126}
\int_{\Omega}\aleph(y,\nabla T_{\tau}(u_{n})).\nabla(T_{\tau}(u_{n})-T_{\tau}(u_{p}))dy
&=\int_{\Omega}(\aleph(y,\nabla T_{\tau}(u_{n}))-\aleph(y,\nabla T_{\tau}(u_{p}))).\nabla(T_{\tau}(u_{n})-T_{\tau}(u_{p}))dy\nonumber\\
&+\int_{\Omega}\aleph(y,\nabla T_{\tau}(u_{p})).\nabla(T_{\tau}(u_{n})-T_{\tau}(u_{p}))dy.
\end{align}
Since $ \aleph(y,\nabla T_{\tau}(u_{p}))\in L^{p^{\prime}}(\Omega)$ and $ T_{\tau}(u_{n}) $ converges to $ T_{\tau}(u_{p}) $ weakly in $W^{1,p}(\Omega),$ we find
\begin{align}\label{79}
&\lim_{n\rightarrow +\infty}\int_{\Omega}\aleph(y,\nabla T_{\tau}(u_{p})).\nabla(T_{\tau}(u_{n})-T_{\tau}(u_{p}))dy=0.
\end{align}
From \eqref{126} and  \eqref{79}, equation \eqref{78} becomes
\begin{equation}\label{81}
\int_{\Omega}\aleph(y,\nabla u_{n}).\nabla(T_{\tau}(u_{n})-T_{\tau}(u_{p}))dy =\int_{\Omega}(\aleph(y,\nabla T_{\tau}(u_{n}))-\aleph(y,\nabla T_{\tau}(u_{p}))).\nabla(T_{\tau}(u_{n})-T_{\tau}(u_{p}))dy+\varepsilon_{3}(n),
\end{equation}
with $\varepsilon_{3}(n)$ is real number such that $\lim\limits_{n\rightarrow +\infty}\varepsilon_{3}(n)=0.$ It follows from  \eqref{127} we obtain that
\begin{align}\label{82}
&\lim_{n\rightarrow +\infty}\int_{\partial \Omega}\lambda \vert u_{n}\vert^{p-2}u_{n}(T_{\tau}(u_{n})-T_{\tau}(u_{p}))d\mathcal{H}^{m-1}=0,
\end{align}
and that
\begin{align}\label{83}
&\lim_{n\rightarrow +\infty}\int_{\partial \Omega}\vert u_{n}\vert^{s}(T_{\tau}(u_{n})-T_{\tau}(u_{p}))d\mathcal{H}^{m-1}=0.
\end{align}
By \eqref{81}, \eqref{82} and \eqref{83}, equation \eqref{77} becomes
\begin{align}\label{128}
&\int_{\Omega}(\aleph(y,\nabla T_{\tau}(u_{n}))-\aleph(y,\nabla T_{\tau}(u_{p}))).\nabla(T_{\tau}(u_{n})-T_{\tau}(u_{p}))dy+\varepsilon_{4}(n)\nonumber \\
&=\int_{\Omega}\frac{f_{n}}{(u_{n}+\frac{1}{n})^{\gamma}}(T_{\tau}(u_{n})-T_{\tau}(u_{p}))dy
+\int_{\partial\Omega}\frac{g_{n}}{(u_{n}+\frac{1}{n})^{\eta}}(T_{\tau}(u_{n})-T_{\tau}(u_{p}))d\mathcal{H}^{m-1},
\end{align}
where $\varepsilon_{4}(n)$ is real number satisfying $\lim\limits_{n\rightarrow +\infty}\varepsilon_{4}(n)=0.$ Regarding the initial term on the right hand side of \eqref{128} we reach that
\begin{equation}\label{64}
\int_{\Omega}\frac{f_{n}}{(u_{n}+\frac{1}{n})^{\gamma}}(T_{\tau}(u_{n})-T_{\tau}(u_{p}))dy\leq\int_{\{u_{n}\leq\delta\}}\frac{f_{n}}{u_{n}^{\gamma}}(T_{\tau}(u_{n})-T_{\tau}(u_{p}))dy+\frac{1}{\delta^{\gamma}}
\int_{\{u_{n}>\delta\}}f_{n}(T_{\tau}(u_{n})-T_{\tau}(u_{p}))dy.
\end{equation}
For the second term on the right hand of \eqref{64} we have that $ f_{n}$ converge  strongly to $ f$ in $ L^{1}(\Omega)$ and $(T_{\tau}(u_{n})-T_{\tau}(u_{p})) $ converges  $ ^{\ast}$- weakly to $ 0$ in $L^{\infty}(\Omega)$  as $ n\rightarrow \infty.$ Hence
\begin{align}\label{75}
\lim_{n\rightarrow +\infty}\frac{1}{\delta^{\gamma}}\int_{\{u_{n}>\delta\}}f_{n}(T_{\tau}(u_{n})-T_{\tau}(u_{p}))dy=0.
\end{align}
For the first term  on the right hand of \eqref{64} we see
\begin{align*}
& \int_{\{u_{n}\leq\delta\}}\frac{f_{n}}{u_{n}^{\gamma}}(T_{\tau}(u_{n})-T_{\tau}(u_{p}))dy\leq \int_{\{u_{n}\leq\delta\}}\delta^{1-\gamma}fdy.
\end{align*}
Since
\begin{align*}
&\chi_{\{u_{n}\leq\delta\}}\delta^{1-\gamma}f\leq \delta^{1-\gamma}f.
\end{align*}
By applying Lebesgue's Theorem, we obtain
\begin{align*}
&\lim_{n\rightarrow +\infty}\int_{\{u_{n}\leq\delta\}}\delta^{1-\gamma}fdy=\int_{\{u_{p}\leq\delta\}}\delta^{1-\gamma}fdy \leq\int_{\{u_{p}\leq\delta\}}fdy.
\end{align*}
It follows from Lemma \ref{99} that $\{u_{p}=0\}\subset \{f=0\} $ up to a set of zero Lebesgue measure, then
\begin{equation}\label{74}
    \lim_{\delta\rightarrow 0^{+}}\int_{\{u_{p}\leq \delta\}}fdy=\int_{\{u_{p}=0\}}fdy=0.
\end{equation}
From \eqref{64}, \eqref{75} and  \eqref{74} we deduce that
\begin{align}\label{129}
&\lim_{\delta\rightarrow 0^{+}}\lim_{n\rightarrow+\infty}\int_{\{u_{n}\leq\delta\}}\frac{f_{n}}{u_{n}^{\gamma}}(T_{\tau}(u_{n})-T_{\tau}(u_{p}))dy\leq 0.
\end{align}
On the other hand we have
\begin{align*}
&\frac{g_{n}}{(u_{n}+\frac{1}{n})^{\eta}}(T_{\tau}(u_{n})-T_{\tau}(u_{p}))\leq 2\tau\frac{g}{\beta_{p}^{\eta}},
\end{align*}
$\mathcal{H}^{m-1}$ a. e.  on $\partial \Omega$. Using Lebesgue's Theorem so that
\begin{align}\label{84}
&\lim_{n\rightarrow+\infty}\int_{\partial \Omega}\frac{g_{n}}{(u_{n}+\frac{1}{n})^{\eta}}(T_{\tau}(u_{n})-T_{\tau}(u_{p}))d\mathcal{H}^{m-1}=0.
\end{align}
By \eqref{d}, \eqref{128}, \eqref{129} and \eqref{84} we arrive at
\begin{equation*}
\lim_{n\rightarrow +\infty}\int_{\Omega}(\aleph(y,\nabla T_{\tau}(u_{n}))-\aleph(y,\nabla T_{\tau}(u_{p})))\nabla(T_{\tau}(u_{n})-T_{\tau}(u_{p}))dy=0.
\end{equation*}
Then \cite[Lemma 5]{6014} implies
\begin{equation*}
    T_{\tau}(u_{n})\rightarrow T_{\tau}(u_{p}) \hspace*{0.2cm} strongly \hspace*{0.1cm} in \hspace*{0.1cm} W^{1,p}(\Omega).
\end{equation*}
\begin{flushright}
    $\square$
\end{flushright}
\end{dem}
\subsection{Proof of Theorem \ref{48}}
\begin{dem}
We use $\varphi \in W^{1,p}(\Omega)\cap L^{\infty}(\Omega)$ as a \textit{t. f.} in \eqref{1} to get
\begin{align}\label{31}
&\int_{\Omega}\aleph(y,\nabla u_{n})\nabla \varphi dy+\int_{\partial \Omega} \lambda \vert u_{n}\vert^{p-2}u_{n}\varphi d\mathcal{H}^{m-1}+\int_{\partial \Omega}\vert u_{n}\vert^{s-1}u_{n}\varphi d\mathcal{H}^{m-1}\nonumber \\
&=\int_{\Omega}\frac{f_{n}}{(u_{n}+\frac{1}{n})^{\gamma}}\varphi dy+\int_{\partial \Omega}\frac{g_{n}}{(u_{n}+\frac{1}{n})^{\eta}}\varphi d\mathcal{H}^{m-1}.
\end{align}
Regarding the first term on the left-hand side of \eqref{31}. From Lemma \textnormal{\ref{96}}, we have, up to subsequences $\nabla u_{n}\rightarrow \nabla u_{p}$ a. e. in $ \Omega$. From this fact and \eqref{c}, which implies
\begin{equation}\label{36}
 \lim_{n\rightarrow +\infty}\int_{\Omega}\aleph(y,\nabla u_{n})\nabla\varphi dy=\int_{\Omega}\aleph(y,\nabla u_{p})\nabla\varphi dy.
\end{equation}
Consider $\delta>0 | \delta\notin \{t:|\{u_{p}=t\}|>0\},$ we can write
\begin{align}\label{10}
&\int_{\Omega}\frac{f_{n}}{(u_{n}+\frac{1}{n})^{\gamma}}\varphi dy=\int_{\{u_{n}\leq \delta\}}\frac{f_{n}}{(u_{n}+\frac{1}{n})^{\gamma}}\varphi dy+\int_{\{u_{n}> \delta\}}\frac{f_{n}}{(u_{n}+\frac{1}{n})^{\gamma}}\varphi dy.
\end{align}
Turning to  the second integral on the right hand side of \eqref{10},  we see that
\begin{align*}
&\chi_{\{u_{n} >\delta\}} \frac{f_{n}}{(u_{n}+\frac{1}{n})^{\gamma}} \leq \frac{f}{\delta^{\gamma}}.
\end{align*}
By Lebesgue's Theorem we have
\begin{align*}
&\lim_{n\rightarrow +\infty}\int_{\{u_{n} >\delta\}} \frac{f_{n}}{(u_{n}+\frac{1}{n})^{\gamma}}\varphi dy
=\int_{\{u_{p} >\delta\}} \frac{f}{u_{p}^{\gamma}}\varphi dy.
\end{align*}
Moreover
\begin{align*}
& \chi_{\{u_{p} >\delta\}}\frac{f}{u_{p}^{\gamma}}\varphi\leq \frac{f}{u_{p}^{\gamma}}\varphi\in L^{1}(\Omega),
\end{align*}
applying Lebesgue's Theorem we get
\begin{align}\label{11}
& \lim_{\delta\rightarrow 0}\lim_{n\rightarrow +\infty}\int_{\{u_{p} >\delta\}} \frac{f}{u_{p}^{\gamma}}\varphi dy=\int_{\{u_{p} >0\}} \frac{f}{u_{p}^{\gamma}}\varphi dy.
\end{align}
Taking $ \mathfrak{S}_{\delta}(u_{n})\varphi $ as a \textit{t. f.} in \eqref{1} we obtain that
\begin{equation*}
\int_{\{u_{n}\leq \delta\}} \frac{f_{n}}{(u_{n}+\frac{1}{n})^{\gamma}}\varphi dy\leq \int_{\Omega}\aleph(y,\nabla u_{n})\mathfrak{S}_{\delta}(u_{n})\nabla \varphi dy+\int_{\partial \Omega}\lambda \vert u_{n}\vert^{p-1} \mathfrak{S}_{\delta}(u_{n})\varphi d\mathcal{H}^{m-1}+\int_{\partial \Omega}\vert u_{n}\vert^{s}\mathfrak{S}_{\delta}(u_{n})\varphi d\mathcal{H}^{m-1}.
\end{equation*}
When $ n\rightarrow +\infty,$ using weak convergence and Lebesgue Theorem,  we reach that
\begin{equation*}
\limsup_{n\rightarrow +\infty}\int_{\{u_{n}\leq \delta\}} \frac{f_{n}}{(u_{n}+\frac{1}{n})^{\gamma}}\varphi dy \leq \int_{\Omega}\aleph(y,\nabla u_{p})\mathfrak{S}_{\delta}(u_{p})\nabla \varphi dy+\int_{\partial \Omega}\lambda \vert u_{p}\vert^{p-1} \mathfrak{S}_{\delta}(u_{p})\varphi d\mathcal{H}^{m-1}+\int_{\partial \Omega}\vert u_{p}\vert^{s}\mathfrak{S}_{\delta}(u_{p})\varphi d\mathcal{H}^{m-1},
\end{equation*}
passing to the limit as $\delta\rightarrow 0^{+},$ it follows that
\begin{align*}
&\lim_{\delta\rightarrow 0^{+} }\limsup_{n\rightarrow +\infty}\int_{\{u_{n}\leq \delta\}} \frac{f_{n}}{(u_{n}+\frac{1}{n})^{\gamma}}\varphi dy\leq \int_{\{u_{p}=0\}}\aleph(y,\nabla u_{p})\mathfrak{S}_{\delta}(u_{p})\nabla \varphi d\mathcal{H}^{m-1}\\
&+\int_{\{u_{p}=0\}}\lambda \vert u_{p}\vert^{p-1} \mathfrak{S}_{\delta}(u_{p})\varphi d\mathcal{H}^{m-1}+\int_{\{u_{p}=0\}}\vert u_{p}\vert^{s}\mathfrak{S}_{\delta}(u_{p})\varphi d\mathcal{H}^{m-1}=0,
\end{align*}
we deduce that
\begin{align}\label{12}
&\lim_{n\rightarrow +\infty}\int_{\{u_{n}\leq \delta\}} \frac{f_{n}}{(u_{n}+\frac{1}{n})^{\gamma}}\varphi dy=0.
\end{align}
By \eqref{11} and \eqref{12} we have
\begin{align}\label{13}
&\lim_{n\rightarrow +\infty}\int_{\Omega} \frac{f_{n}}{(u_{n}+\frac{1}{n})^{\gamma}}\varphi dy=\int_{\{u_{p}>0\}} \frac{f}{u_{p}^{\gamma}}\varphi dy=\int_{\Omega} \frac{f}{u_{p}^{\gamma}}\varphi dy.
\end{align}
Since
\begin{align*}
&\frac{g_{n}}{(u_{n}+\frac{1}{n})^{\eta}}\leq \frac{g}{\beta_{p}^{\eta}},
\end{align*}
utilizing Lebesgue's Theorem we arrive at
\begin{align}\label{14}
& \lim_{n\rightarrow +\infty}\int_{\partial \Omega}\frac{g_{n}}{(u_{n}+\frac{1}{n})^{\eta}}d\mathcal{H}^{m-1}=\int_{\Omega}\frac{g}{u_{p}^{\eta}}d\mathcal{H}^{m-1}.
\end{align}
Letting $n$ tends to $+\infty $ in \eqref{31}, it follows from  \eqref{38}, \eqref{36}, \eqref{13} and \eqref{14} we conclude that \eqref{19}.
\begin{flushright}
    $\square$
\end{flushright}
\end{dem}
\section{Existence and uniqueness of solution for $ p=1$ and $0<f\in L^{m}(\Omega)$}\label{sec5}
This section is devoted to solve problem
\begin{equation}\label{4}
    \left\{
  \begin{array}{ll}
   -\Delta_{1}u=\frac{f}{u^{\gamma}}& \hbox{in $\Omega,$} \\
  \frac{Du}{\vert Du\vert}.\sigma+\frac{u}{\vert u\vert}+\vert u\vert^{s-1}u=0 & \hbox{on $\partial\Omega,$} \\
  \end{array}
\right.
\end{equation}
where $ f$ is positive function in $ L^{m}(\Omega),$ $ 0<\gamma\leq 1$ and $s\geq 1.$
\subsection{Existence result}
We give the following  notion of  weak solution to problem \eqref{4}.
\begin{defn}
Let $0<f\in L^{m}(\Omega).$ A function $ u\in BV(\Omega)$ is  a weak solution to the problem \eqref{4}, if  there exist  vector field  $z\in \mathcal{DM}^{\infty}_{loc}(\Omega),$  $ v\in L^{s+1}(\partial \Omega)$ and $\vartheta\in L^{\infty}(\partial\Omega)$  such that
\begin{eqnarray}
& & \frac{f}{u^{\gamma}}\in L^{1}_{loc}(\Omega), \label{22} \\
& & \vert \vert z\vert \vert _{L^{\infty}(\Omega)^{m}}\leq 1,\\
& & -\textnormal{div }z=\frac{f}{u^{\gamma}}\hspace*{1mm} \mbox{in} \hspace*{1mm} \mathcal{D}^{\prime}(\Omega),\label{23} \\
& & (z,DT_{\tau}(u))=\vert DT_{\tau}(u)\vert   \hspace*{1mm} \mbox{as}\hspace*{1mm} \mbox{measures} \hspace*{1mm}\mbox{in}\hspace*{1mm} \Omega \mbox{ for any } \tau>0,\label{24}\\
& &  T_{\tau}(u)+T_{\tau}(v)^{s+1}+T_{\tau}(u)[z,\sigma]=0\hspace*{0.4cm}  \mathcal{H}^{m-1}- \mbox{a. e. on } \partial \Omega\mbox{ for any } \tau>0,\label{25}\\
& & \vartheta + v^{s}+[z,\sigma]=0 \hspace*{0.4cm}\mathcal{H}^{m-1}-\mbox{ a. e. on } \partial \Omega.\label{44}
\end{eqnarray}
\end{defn}
Now we can state our theorem:
\begin{thm}\label{93}
Assume that $ 0<\gamma\leq 1,$ $s\geq 1$ and $0<f\in L^{m}(\Omega).$  Then there exists a solution to \eqref{4}.
\end{thm}
In order to prove Theorem \ref{93}, we will use a sequence  of Lemmas. Let  $ 1<p<2,$ we need the following approximation:
\begin{equation}\label{55}
    \left\{
  \begin{array}{ll}
   -\Delta_{p}u_{p}=\frac{f}{u_{p}^{\gamma}}& \hbox{in $\Omega,$} \\
  \frac{\partial u_{p}}{\partial \sigma}+\vert u_{p}\vert^{p-2} u_{p}+\vert u_{p}\vert^{s-1}u_{p}=0 & \hbox{on $\partial\Omega.$} \\
  \end{array}
\right.
\end{equation}
Thanks to Theorem \ref{48},  problem \eqref{55} has a solution $ u_{p}\in W^{1,p}(\Omega).$\\

Under the previous assumptions, we have a constant  $ \beta_{p}$  that depends on $ p$ and may possibly  vanish as $ p$ tends to 1. Moreover, we will lack pointwise convergence of $ u_{p}$ on $ \partial \Omega$ as $ p$ approaches 1. For these reasons, we set $ g=0$ in \eqref{53}.\\

In this section, we are interested in the behavior as $ p$ tends to $ 1$ in \eqref{55}.\\

\begin{lem}\label{5}
Let $ f\in L^{m}(\Omega)$ be nonnegative, $ s\geq 1$  and $ 0<\gamma \leq1,$  let $ u_{p}$ be a solution of \eqref{55}. Then there exist   functions $ u\in BV(\Omega),$ $v\in L^{s+1}(\partial\Omega) $ and $\vartheta\in L^{\infty}(\partial\Omega)$ such that, up to subsequences, satisfying
\begin{eqnarray}
& & u_{p} \rightarrow u \hspace*{2mm}\mbox{strongly} \hspace*{1mm} \mbox{in} \hspace*{1mm} L^{b}(\Omega) \hspace*{3mm}1\leq b<\frac{m}{m-1} \hspace*{1mm} \mbox{as} \hspace*{1mm} p\rightarrow 1^{+},\label{88} \\
& & \nabla u_{p}\overset{\ast}{\rightharpoonup} Du \hspace*{1mm}\mbox{as} \hspace*{1mm} \mbox{measures} \hspace*{1mm} \mbox{as} \hspace*{1mm}p\rightarrow 1^{+},  \label{89}\\
& &  u\geq 0 \hspace*{2mm}\mbox {a. e. in }\hspace*{1mm} \Omega,\label{90}\\
& & u_{p}\rightharpoonup v \hspace*{2mm} weakly \hspace*{2mm} in \hspace*{2mm} L^{s+1}(\partial \Omega)\hspace*{1mm}\mbox{as} \hspace*{1mm}p\rightarrow 1^{+},\label{91}\\
& &  v\geq 0 \hspace*{2mm}\mbox { a. e. on }\hspace*{1mm} \partial\Omega,\label{92}\\
& & \vert u_{p}\vert^{p-1}u_{p}\rightharpoonup \vartheta \hspace*{2mm} weakly \hspace*{2mm} in \hspace*{2mm} L^{\varrho}(\partial \Omega),\hspace*{1mm} \forall \varrho> 1\hspace*{1mm}\mbox{as} \hspace*{1mm}p\rightarrow 1^{+}.
\end{eqnarray}
Moreover,
\begin{equation*}
 \Vert \vartheta \Vert_{L^{\infty}(\partial \Omega)}\leq 1.
\end{equation*}
\end{lem}
\begin{dem}
Using lower semicontinuity of the norm concerning  $n$  in \eqref{40}, we get
\begin{equation}\label{122}
    \Vert u_{p} \Vert_{W^{1,p}(\Omega)} \leq C_{13},
\end{equation}
with $ C_{13}$ is  independent of $p.$  With  the assistance  of  Young's inequality , it seems that
\begin{equation}\label{123}
  \Vert u_{p}\Vert_{W^{1,1}(\Omega)}=  \int_{\Omega}\vert \nabla u_{p}\vert dy+\int_{\Omega}\vert u_{p}\vert dy\leq \frac{1}{p}\int_{\Omega}\vert \nabla u_{p}\vert^{p} dy+\frac{1}{p}\int_{\Omega}\vert u_{p}\vert^{p} dy+2\frac{p-1}{p}\vert \Omega \vert\leq C_{14}.
\end{equation}
where $ C_{14}$ does not  depend  of $ p.$  As a result, we are able to identify a function $ u\in BV(\Omega)$ satisfying, up to subsequence, \eqref{88}, \eqref{89} and \eqref{90}.
From  lower semicontinuity of the norm concerning $ n$ in \eqref{43} we obtain
 \begin{equation*}
    \int_{\partial \Omega}\vert u_{p}\vert^{s+1} d\mathcal{H}^{m-1}\leq C_{15},
 \end{equation*}
with $ C_{15}$ is  independent of $ p.$  Therefore, there exists $ v\in L^{s+1}(\partial \Omega)$ satisfying, up to subsequence,
\begin{equation*}
 u_{p}\rightharpoonup v \hspace*{2mm} weakly \hspace*{2mm} in \hspace*{2mm} L^{s+1}(\partial \Omega)\hspace*{1mm}\mbox{as} \hspace*{1mm}p\rightarrow 1^{+} \,\mbox{ and } \, v\geq 0 \hspace*{2mm}\mbox {a. e. on}\hspace*{1mm}\partial\Omega.
\end{equation*}
 Let $ \varrho>1,$ we can consider  $1<p<\frac{\varrho}{\varrho-1},$ by H\"{o}lder's inequality, it follows that
\begin{eqnarray}\label{100}
\Vert u_{p}^{p-1} \Vert_{L^{\varrho}(\partial \Omega)}\leq \bigg( \int_{\partial \Omega}\vert u_{p} \vert ^{p} d\mathcal{H}^{m-1}\bigg)^{\frac{p-1}{p}}  \bigg(\mathcal{H}^{m-1}(\partial \Omega)\bigg)^{\frac{1}{\varrho}-\frac{p-1}{p}}\leq C_{16},
\end{eqnarray}
where $C_{16}$ does not depend on $p.$ Thus, up to subsequence there exists $\vartheta_{\varrho}\in L^{\varrho}(\partial \Omega)$ such that
\begin{equation*}
    \vert u_{p}\vert ^{p-2}u_{p}\rightharpoonup  \vartheta_{\varrho}\hspace*{2mm} weakly \hspace*{2mm} in \hspace*{2mm} L^{\varrho}(\partial \Omega)\hspace*{1mm}\mbox{as} \hspace*{1mm}p\rightarrow 1^{+}.
\end{equation*}
It follows from a diagonal argument that there exists $ \vartheta$ satisfying
\begin{equation*}
    \vert u_{p}\vert ^{p-2}u_{p}\rightharpoonup  \vartheta\hspace*{2mm} weakly \hspace*{2mm} in \hspace*{2mm} L^{\varrho}(\partial \Omega)\hspace*{1mm}\mbox{as} \hspace*{1mm}p\rightarrow 1^{+}.
\end{equation*}
Using the lower semicontinuity of the $\varrho$-norm  in \eqref{100}, we get
\begin{equation}\label{101}
    \Vert \vartheta \Vert_{L^{\varrho}(\partial \Omega)}\leq \bigg(\mathcal{H}^{m-1}(\partial \Omega)\bigg)^{\frac{1}{\varrho}}.
\end{equation}
Taking $ \varrho$ tends to $+\infty$ in \eqref{101}, it yields
\begin{equation}\label{102}
\Vert \vartheta \Vert_{L^{\infty}(\partial \Omega)}\leq 1.
\end{equation}
\begin{flushright}
    $\square$
\end{flushright}
\end{dem}
\begin{lem}\label{6}
Assume that $f$ belongs to $L^{m}(\Omega)$ is positive, \( s \geq 1 \), and \( 0 < \gamma \leq 1 \). Let \( u \) be the function obtained in Lemma \ref{5}. Then  there exists a vector field $z\in \mathcal{DM}^{\infty}_{loc}(\Omega)$ satisfying
\begin{eqnarray}
   & & \frac{f}{u^{\gamma}}\in L^{1}_{loc}(\Omega), \nonumber \\ 
   & & \Vert z\Vert_{L^{\infty}(\Omega)^{m}}\leq 1, \nonumber\\
   & & -\textnormal{div }z=\frac{f}{u^{\gamma}}\hspace*{1mm} \mbox{in} \hspace*{1mm} \mathcal{D}^{\prime}(\Omega),\nonumber\\ 
   & & (z,DT_{\tau}(u))=\vert DT_{\tau}(u)\vert   \hspace*{1mm} \mbox{as}\hspace*{1mm} \mbox{measures} \hspace*{1mm}\mbox{in}\hspace*{1mm} \Omega \mbox{ for any } \tau>0.\nonumber 
\end{eqnarray}
\end{lem}
\begin{dem}
\textnormal{\textbf{Step 1:}} $\frac{f}{u^{\gamma}}\in L^{1}_{loc}(\Omega).$\\
Employing  Fatou's Lemma in \eqref{41} leads to
\begin{align*}
& \int_{\Omega}\frac{f}{u^{\gamma}}\varphi dy\leq \liminf_{p\rightarrow 1^{+}}\int_{\Omega}\frac{f}{u_{p}^{\gamma}}\varphi dy\leq C_{17},
\end{align*}
Hence
\begin{equation*}
    \frac{f}{u^{\gamma}}\in L^{1}_{loc}(\Omega).
\end{equation*}
Consequently we conclude that, up to a set of zero Lebesgue measure
\begin{equation*}
    \{u=0\}\subset \{f=0\}.
\end{equation*}
\textnormal{\textbf{Step 2:}} Existence $ z\in\mathcal{DM}^{\infty}_{loc}(\Omega)$ with $\Vert z\Vert_{L^{\infty}(\Omega)^{m}}\leq 1$ and
$ -\textnormal{div }z=\frac{f}{u^{\gamma}}$  holds in $ \mathcal{D}^{\prime}(\Omega).$
\par We will demonstrate that the sequence  $\{\vert \nabla u_{p}\vert ^{p-2}\nabla u_{p}\}_{p>1}$ is bounded in  $L^{\xi}(\Omega;\mathbb{R}^{m})$ with $ \xi\geq 1.$ To achieve this, fix $ 1\leq \xi<+\infty$ and take $ 1<p<\frac{\xi}{\xi-1}.$ Utilizing H$\ddot{o}$lder's inequality and \eqref{122}, our conclusion is that
\begin{equation}\label{201}
    \vert \vert \vert \nabla u_{p}\vert ^{p-2}\nabla u_{p}\vert \vert _{L^{\xi}(\Omega;\mathbb{R}^{m})}\leq \left(\int_{\Omega}\vert \nabla u_{p}\vert ^{p}\right)^{\frac{1}{p^{\prime}}}\vert \Omega\vert ^{\frac{1}{\xi}-\frac{1}{p^{\prime}}}\leq (C_{13})^{p-1}\vert \Omega\vert ^{\frac{1}{\xi}-\frac{1}{p^{\prime}}}.
\end{equation}
Hence, from \eqref{201}  we find $ z\in L^{\xi}(\Omega;\mathbb{R}^{m})$ where  
\begin{equation}\label{200}
   \vert \nabla u_{p}\vert ^{p-2}\nabla u_{p} \mbox{ converges weakly to } z \mbox{ in } L^{\xi}(\Omega;\mathbb{R}^{m}) \mbox{ as } p\rightarrow 1^{+} \mbox{ and }  \xi\geq 1.
\end{equation}
Moreover letting $p$ as $ 1^{+}$ in \eqref{201}, using  the lower semicontinuity it produces
\begin{equation*}
    \vert \vert z\vert \vert _{L^{\xi}(\Omega;\mathbb{R}^{m})}\leq \vert \Omega\vert ^{\frac{1}{\xi}}, \hspace*{0.2cm}  \forall \xi<+\infty,
\end{equation*}
passing to the limit as  $ \xi\rightarrow+\infty,$  we have  $ \Vert z \Vert _{L^{\infty}(\Omega)^{m}}\leq 1.$
\par Using $0 \leq \phi \in C^{\infty}_{0}(\Omega)$ as a \textit{t. f.} in \eqref{55}, we obtain
\begin{align}\label{21}
&\int_{\Omega}\vert \nabla u_{p}\vert^{p-2}\nabla u_{p}.\nabla \phi dy=\int_{\Omega}\frac{f}{u^{\gamma}_{p}}\phi dy.
\end{align}
Take $\delta>0 $ with   $ \delta \notin \{t:\vert \{u_{p}=t\}\vert >0\},$ we decompose the term in equation \eqref{21} as follows
\begin{align}\label{104}
&\int_{\Omega}\frac{f}{u^{\gamma}_{p}}\phi dy=\int_{\{u_{p}\leq \delta\}}\frac{f}{u^{\gamma}_{p}}\phi dy+\int_{\{u_{p}> \delta\}}\frac{f}{u^{\gamma}_{p}}\phi dy.
\end{align}
Notice  that
\begin{align*}
&\chi_{\{u_{p}> \delta\}}\frac{f}{u^{\gamma}_{p}}\leq\frac{f}{\delta^{\gamma}}.
\end{align*}
Using Lebesgue's Theorem we have
\begin{align*}
&\lim_{p\rightarrow 1^{+}}\int_{\{u_{p}> \delta\}}\frac{f}{u^{\gamma}_{p}}\phi dy=\int_{\{u> \delta\}}\frac{f}{u^{\gamma}}\phi dy.
\end{align*}
We remark that
\begin{equation*}
  \chi_{\{u> \delta\}} \frac{f}{u^{\gamma}}\leq\frac{f}{u^{\gamma}}.
\end{equation*}
Once more applying Lebesgue's Theorem we obtain
\begin{align}\label{27}
&\lim_{\delta\rightarrow 0^{+}}\int_{\{u> \delta\}}\frac{f}{u^{\gamma}}\phi dy=\int_{\{u>0\}}\frac{f}{u^{\gamma}}\phi dy.
\end{align}
To handle the first term on the right-hand side of \eqref{104},  taking $\mathfrak{S}_{\delta}(u_{p})\phi $  as a \textit{t. f.} in \eqref{55} with $0\leq\phi\in C^{\infty}_{0}(\Omega),$ we get
\begin{align*}
\int_{\{u_{p}\leq \delta\}}\frac{f}{u^{\gamma}_{p}}\phi dy &\leq\int_{\Omega}\frac{f}{u^{\gamma}_{p}}\mathfrak{S}_{\delta}(u_{p})\phi dy=\int_{\Omega}\vert \nabla u_{p}\vert ^{p-2}\nabla u_{p}\mathfrak{S}_{\delta}(u_{p})\nabla \phi dy-\int_{\Omega}\frac{1}{\delta}\vert \nabla u_{p}\vert ^{p}\phi dy\\
&\leq \int_{\Omega}\vert \nabla u_{p}\vert ^{p-2}\nabla u_{p}\mathfrak{S}_{\delta}(u_{p})\nabla \phi dy.
\end{align*}
Thanks to \eqref{200}  we arrive at
\begin{align*}
\limsup_{p\rightarrow 1^{+}}\int_{\{u_{p}\leq \delta\}}\frac{f}{u^{\gamma}_{p}}\phi dy \leq \int_{\Omega}z\mathfrak{S}_{\delta}(u)\nabla \phi dy.
\end{align*}
We take the limit in $ \delta\rightarrow 0^{+},$ so that
\begin{align}\label{26}
\lim_{\delta\rightarrow 0^{+}}\limsup_{p\rightarrow 1^{+}}\int_{\{u_{p}\leq \delta\}}\frac{f}{u^{\gamma}_{p}}\phi dy \leq \int_{\{u=0\}}z\nabla \phi dy=0.
\end{align}
From \eqref{27} and  \eqref{26}  we have
\begin{align}\label{105}
\lim_{p\rightarrow 1^{+}}\int_{\Omega}\frac{f}{u^{\gamma}_{p}}\phi dy=\int_{\{u>0\}}\frac{f}{u^{\gamma}} \phi dy=\int_{\Omega}\frac{f}{u^{\gamma}} \phi dy.
\end{align}
Letting $ p \rightarrow 1^{+}$ in \eqref{21}, by \eqref{200} and \eqref{105}  we reach at
\begin{align*}
\int_{\Omega}z \nabla \phi dy=\int_{\Omega}\frac{f}{u^{\gamma}}\phi dy.
\end{align*}
As a result \eqref{23}.
\par From \eqref{23} we conclude that $ z\in\mathcal{DM}^{\infty}_{loc}(\Omega).$\\
\textnormal{\textbf{Step 3:}} $(z,DT_{\tau}(u))=\vert DT_{\tau}(u)\vert$ as measures in $ \Omega.$
\par Choosing $(T_{\tau}(u)\ast \rho_{\epsilon})\varphi$ as a \textit{t. f.} in \eqref{23} where $\rho_{\epsilon} $ is a standard mollifier and $\varphi$ belongs to $C^{\infty}_{0}(\Omega).$ Thus
\begin{equation}\label{42}
  -\int_{\Omega}(T_{\tau}(u)\ast \rho_{\epsilon})\varphi \textnormal{div } z= \int_{\Omega}\frac{f}{u^{\gamma}}(T_{\tau}(u)\ast \rho_{\epsilon})\varphi dy.
\end{equation}
Since $T_{\tau}(u)\ast \rho_{\epsilon}\rightarrow T_{\tau}(u) $ $\mathcal{H}^{m-1}-$ a. e. as $\varepsilon\rightarrow 0^{+}$ and
$T_{\tau}(u)\ast \rho_{\epsilon}\leq \tau.$ We obtain by letting   $\varepsilon\rightarrow 0^{+}$ in \eqref{42} that
\begin{equation*}
   -\int_{\Omega}T_{\tau}(u)^{\ast}\varphi\textnormal{div } z= \int_{\Omega}\frac{f}{u^{\gamma}}T_{\tau}(u)^{\ast}\varphi dy.
\end{equation*}
Hence
\begin{equation*}
 - T_{\tau}(u)^{\ast}\textnormal{div } z= \frac{f}{u^{\gamma}}T_{\tau}(u)^{\ast} \hspace*{0.3cm}in \hspace*{0.3cm}\mathcal{D}^{\prime}(\Omega).
\end{equation*}
Let $\tau>0, $  $ 0\leq\varphi \in C^{\infty}_{0}(\Omega),$  choosing $ T_{\tau}(u_{p})\varphi$ as a  \textit{t. f.} in \eqref{55}, Young's inequality imply that
\begin{equation}\label{108}
\int_{\Omega}\vert \nabla T_{\tau}(u_{p})\vert \varphi dy+\int_{\Omega}T_{\tau}(u_{p})\vert \nabla u_{p}\vert^{p-2}\nabla u_{p}\nabla \varphi dy \leq\int_{\Omega} \frac{f}{u^{\gamma}_{p}}T_{\tau}(u_{p})\varphi dy+\frac{p-1}{p}\int_{\Omega}\varphi dy.
\end{equation}
If $ 0<\gamma\leq1,$   we see that
\begin{equation*}
\frac{f}{u_{p}^{\gamma}}T_{\tau}(u_{p})\leq \tau^{1-\gamma}f.
\end{equation*}
Consequently, by  Lebesgue's Theorem we have
\begin{equation}\label{54}
  \lim_{p\rightarrow 1^{+}}\int_{\Omega}\frac{f}{u_{p}^{\gamma}}T_{\tau}(u_{p})\varphi dy=\int_{\Omega}\frac{f}{u^{\gamma}}T_{\tau}(u)\varphi dy.
\end{equation}
Expressions \eqref{108} and  \eqref{54} yield
\begin{align*}
\int_{\Omega}\varphi\vert DT_{\tau}(u)\vert +\int_{\Omega}T_{\tau}(u)z\nabla \varphi dy \leq \int_{\Omega}\frac{f}{u^{\gamma}}T_{\tau}(u)\varphi dy=-\int_{\Omega}T_{\tau}(u)^{\ast} \varphi \textnormal{div } z.
\end{align*}
Thanks to \eqref{95}, we obtain
\begin{equation*}
    \int_{\Omega}\varphi \vert DT_{\tau}(u)\vert\leq \int_{\Omega}\varphi (z,DT_{\tau}(u)).
\end{equation*}
Therefore
\begin{equation*}
     \vert DT_{\tau}(u)\vert\leq (z,DT_{\tau}(u))\hspace*{0.1cm}as\hspace*{0.1cm}measures\hspace*{0.1cm}in \hspace*{0.1cm}\Omega.
\end{equation*}
Since $ \Vert z\Vert_{L^{\infty}(\Omega)^{m}}\leq 1 $ the reverse inequality is valid, namely
\begin{equation*}
(z,DT_{\tau}(u))=\vert DT_{\tau}(u)\vert \hspace*{0.1cm}as\hspace*{0.1cm}measures\hspace*{0.1cm}in \hspace*{0.1cm}\Omega \mbox{ for any } \tau>0.
\end{equation*}
\begin{flushright}
    $\square$
\end{flushright}
\end{dem}
\begin{lem}\label{7}
Let $0<f\in L^{m}(\Omega),$ $ 0<\gamma\leq 1$ and  $s\geq 1.$ Suppose $ u$  and $ v$ are the functions in Lemma \ref{5} and  $ z$ is mentioned  in Lemma \ref{6}. Then
\begin{eqnarray*}
& & T_{\tau}(u)+ T_{\tau}(v)^{s+1}+T_{\tau}(u)[z,\sigma]=0\hspace*{0.4cm}\mathcal{H}^{m-1}-\mbox{ a. e. on }\hspace*{1mm} \partial \Omega \mbox{ for any } \tau>0,\\
& & \vartheta+v^{s}+[z,\sigma]=0 \hspace*{0.4cm}\mathcal{H}^{m-1}-\mbox{ a. e. on }\hspace*{1mm} \partial \Omega.
\end{eqnarray*}
\end{lem}
\begin{dem}
\par Let $\tau>0,$ using $T_{\tau}(u_{p})$ as a \textit{t. f.} in \eqref{55} and applying Young's inequality, we get
\begin{align}\label{109}
&\int_{\Omega}\vert \nabla T_{\tau}(u_{p})\vert dy+\int_{\partial \Omega}T_{\tau}(u_{p}) d\mathcal{H}^{m-1}+\int_{\partial \Omega}T_{\tau}(u_{p})^{s+1}d\mathcal{H}^{m-1}\nonumber\\
&\leq\int_{\Omega}\frac{f}{u_{p}^{\gamma}}T_{\tau}(u_{p})dy+\frac{p-1}{p}\int_{\Omega} dy+\frac{p-1}{p}\int_{\partial \Omega} d\mathcal{H}^{m-1}.
\end{align}
To pass to the limit in \eqref{109} as $p\rightarrow 1^{+},$ by utilizing the lower semicontinuity and \eqref{91} on the left-hand side, and applying Lebesgue's Theorem on the right-hand side,  we arrive at
\begin{equation*}
\int_{\Omega}\vert DT_{\tau}(u)\vert+\int_{\partial \Omega}T_{\tau}(u)d\mathcal{H}^{m-1}+\int_{\partial \Omega}T_{\tau}(v)^{s+1}d\mathcal{H}^{m-1}\leq \int_{\Omega}\frac{f}{u^{\gamma}}T_{\tau}(u)dy.
\end{equation*}
By Lemma \ref{106} and \eqref{24}, we have
\begin{equation*}
\int_{\partial \Omega}(T_{\tau}(u)+T_{\tau}(u)[z,\sigma]+T_{\tau}(v)^{s+1})d\mathcal{H}^{m-1}\leq 0.
\end{equation*}
Therefore
\begin{equation}\label{113}
T_{\tau}(u)+ T_{\tau}(v)^{s+1}+T_{\tau}(u)[z,\sigma]\leq 0\hspace*{0.4cm}\mathcal{H}^{m-1}-\mbox{ a. e. }\hspace*{1mm}on \hspace*{1mm} \partial \Omega \mbox{ for any } \tau>0.
\end{equation}
Since  $\vert[T_{\tau}(u)z,\sigma]\vert \leq T_{\tau}(u)$ we reach at
\begin{align}\label{112}
& 0\leq T_{\tau}(u)+ T_{\tau}(v)^{s+1}+T_{\tau}(u)[z,\sigma]\hspace*{0.3cm} on\hspace*{0.2cm}\partial \Omega.
\end{align}
From \eqref{113} and \eqref{112} we arrive at
\begin{equation*}
 T_{\tau}(u)+ T_{\tau}(v)^{s+1}+T_{\tau}(u)[z,\sigma]=0\hspace*{0.4cm}\mathcal{H}^{m-1}-\mbox{ a. e. }\hspace*{1mm}on \hspace*{1mm} \partial \Omega.
\end{equation*}
 For $ p\in]1,2[,$ let $ \varphi \in W^{1,2}(\Omega)\cap L^{\infty}(\Omega)$  as a \textit{t. f.} in \eqref{55} we possess
\begin{equation*}
\int_{\Omega}\vert \nabla u_{p}\vert^{p-2}\nabla u_{p}\nabla \varphi dy+\int_{\partial \Omega}\vert u_{p}\vert^{p-2}u_{p}\varphi d\mathcal{H}^{m-1}+\int_{\partial \Omega} u_{p}^{s}\varphi d\mathcal{H}^{m-1}=\int_{\Omega} \frac{f}{u_{p}^{\gamma}}\varphi dy.
\end{equation*}
When  $ p$ tends to $1^{+},$ it leads to
\begin{align}\label{136}
\int_{\Omega}z\nabla \varphi dy+\int_{\partial \Omega}\vartheta\varphi d\mathcal{H}^{m-1}+\int_{\partial \Omega} v^{s}\varphi d\mathcal{H}^{m-1} =\int_{\Omega} \frac{f}{u^{\gamma}}\varphi dy.
\end{align}
Moreover, \eqref{136} can be generalized for all $ \varphi \in W^{1,1}(\Omega)\cap L^{\infty}(\Omega).$\\
Utilizing \eqref{23} and  Lemma \ref{106} we have
\begin{align*}
\int_{\partial \Omega}\varphi [z,\sigma]d\mathcal{H}^{m-1}+\int_{\partial \Omega}\vartheta\varphi d\mathcal{H}^{m-1}+\int_{\partial \Omega} v^{s}\varphi d\mathcal{H}^{m-1} =0.
\end{align*}
for every $\varphi \in W^{1,1}(\Omega)\cap L^{\infty}(\Omega),$ additionally it is valid for all $ \varphi \in L^{1}(\partial\Omega)\cap L^{\infty}(\Omega).$ Therefore
\begin{align*}
\vartheta+v^{s}+[z,\sigma]=0 \hspace*{0.1cm}\mathcal{H}^{m-1}\hspace*{0.1cm}-\mbox{ a. e. }\hspace*{0.1cm} on\hspace*{0.1cm}\partial \Omega.
\end{align*}
\begin{flushright}
    $\square$
\end{flushright}
\end{dem}
\begin{proof}(Proof of Theorem \ref{93})\\
Let $ u_{p}$ be a solution of \eqref{55}, by   Lemma \ref{5}  there exist  functions $ u\in BV(\Omega),$ $ v\in L^{s+1}(\partial \Omega)$ and $\vartheta\in L^{\infty}(\partial\Omega)$  such that up to subsequence $u_{p}$ converges a. e. to $ u$ in $\Omega $ as $p\rightarrow 1^{+},$  $u_{p}$ converges weakly to $ v$ in $L^{s+1}(\partial \Omega)$ as $p\rightarrow 1^{+}$ and $\vert u_{p}\vert^{p-1}u_{p}\rightharpoonup \vartheta$ in $L^{\varrho}(\partial \Omega),$ $\forall \varrho>1$ as $ p\rightarrow 1^{+}.$ By Lemma \ref{6}, it can be observed that $ z$ is an element of $\mathcal{DM}_{loc}^{\infty}(\Omega)$ where $ \Vert z\Vert_{L^{\infty}(\Omega)^{m}}\leq 1,$ as a further point
the following equalities hold  \eqref{22}, \eqref{23} and \eqref{24}. Finally Lemma \ref{7} gives  \eqref{25} and \eqref{44}.
\end{proof}
\subsection{Uniqueness of the Solution in $ \Omega$}
Our main uniqueness result in this case is the following.
\begin{thm}\label{87}
Assume \( 0 < f \in L^{m}(\Omega) \), \( 0 < \gamma \leq 1 \), and \( s \geq 1 \). If \( u_1 \) and \( u_2 \) are two solutions of problem \eqref{4}, then \( u_1 = u_2 \) a. e. in \( \Omega \).

\end{thm}
\begin{dem}
Suppose that  $ u_{i}$ \textnormal{(i=1,2)} are two solutions of problem \eqref{4} such that  $z_{i}\in\mathcal{DM}^{\infty}_{loc}(\Omega),$
$ \Vert z_{i} \Vert_{L^{\infty}(\Omega)^{m}}\leq 1$ and
\begin{eqnarray}
& & -\textnormal{div }z_{i}=\frac{f}{u_{i}^{\gamma}}\hspace*{1mm} \mbox{in} \hspace*{1mm} \mathcal{D}^{\prime}(\Omega),\label{66} \\
& & (z_{i},DT_{\tau}(u_{i}))=\vert DT_{\tau}(u_{i})\vert   \hspace*{1mm} \mbox{as}\hspace*{1mm} \mbox{measures} \hspace*{1mm}\mbox{in}\hspace*{1mm} \Omega \mbox{ for any } \tau>0,\\ 
& &  T_{\tau}(u_{i})+T_{\tau}(v_{i})^{s+1}+T_{\tau}(u_{i})[z_{i},\sigma]=0\hspace*{0.4cm}  \mathcal{H}^{m-1}- \mbox{ a. e. on } \partial \Omega \mbox{ for any } \tau>0,\label{68}\\
& &\vartheta_{i}+v_{i}^{s}+[z_{i},\sigma]=0 \hspace*{0.4cm}\mathcal{H}^{m-1}-\mbox{ a. e. on } \partial \Omega.\\ 
\end{eqnarray}
 \par Let $\tau>0,$ taking the \textit{t. f.} $(T_{\tau}(u_{1})-T_{\tau}(u_{2}))$ in the difference of  the weak formulation \eqref{66} solved by
$ u_{1}$ and $ u_{2},$ we obtain
\begin{align*}
&\int_{\Omega}\vert DT_{\tau}(u_{1})\vert-\int_{\Omega}(z_{2},DT_{\tau}(u_{1}))+\int_{\Omega}\vert DT_{\tau}(u_{2})\vert-\int_{\Omega}(z_{1},DT_{\tau}(u_{2}))\\
&-\int_{\partial \Omega}[T_{\tau}(u_{1})z_{1},\sigma]d\mathcal{H}^{m-1}-\int_{\partial \Omega}[T_{\tau}(u_{2})z_{2},\sigma]d\mathcal{H}^{m-1}+\int_{\partial \Omega}[T_{\tau}(u_{2})z_{1},\sigma]d\mathcal{H}^{m-1}\\
&+\int_{\partial\Omega}[T_{\tau}(u_{1})z_{2},\sigma]d\mathcal{H}^{m-1}
=\int_{\Omega}\bigg(\frac{f}{u^{\gamma}_{1}}-\frac{f}{u^{\gamma}_{2}}\bigg)(T_{\tau}(u_{1})-T_{\tau}(u_{2}))dy.
\end{align*}
From \eqref{68}, it follows that
\begin{align*}
&\int_{\Omega}\vert DT_{\tau}(u_{1})\vert-\int_{\Omega}(z_{2},DT_{\tau}(u_{1}))+\int_{\Omega}\vert DT_{\tau}(u_{2})\vert-\int_{\Omega}(z_{1},DT_{\tau}(u_{2}))\\
&+\int_{\partial \Omega}\bigg(T_{\tau}(u_{1})+T_{\tau}(v_{1})^{s+1}-T_{\tau}(u_{1})[z_{2},\sigma]\bigg)d\mathcal{H}^{m-1}\\
&+\int_{\partial \Omega}\bigg(T_{\tau}(u_{2})+T_{\tau}(v_{2})^{s+1}-T_{\tau}(u_{2})[z_{1},\sigma]\bigg)d\mathcal{H}^{m-1}\\
&=\int_{\Omega}\bigg(\frac{f}{u^{\gamma}_{1}}-\frac{f}{u^{\gamma}_{2}}\bigg)(T_{\tau}(u_{1})-T_{\tau}(u_{2}))dy
=\int_{\Omega}\frac{1}{(u_{1}u_{2})^{\gamma}}f\bigg(u^{\gamma}_{2}-u^{\gamma}_{1}\bigg)\bigg(T_{\tau}(u_{1})-T_{\tau}(u_{2})\bigg)dy\leq 0.
\end{align*}
Using that $ \vert T_{\tau}(u_{j})[z_{i},\sigma]\vert\leq T_{\tau}(u_{j}) $ for $i,j=1,2$  we find that
\begin{align}\label{70}
&\int_{\partial \Omega}\bigg(T_{\tau}(u_{1})+T_{\tau}(v_{1})^{s+1}-T_{\tau}(u_{1})[z_{2},\sigma]\bigg)d\mathcal{H}^{m-1}\geq0,
\end{align}
and that
\begin{align}\label{71}
&\int_{\partial \Omega}\bigg(T_{\tau}(u_{2})+T_{\tau}(v_{2})^{s+1}-T_{\tau}(u_{2})[z_{1},\sigma]\bigg)d\mathcal{H}^{m-1}\geq 0.
\end{align}
Since
$\displaystyle{\int_{\Omega}}\bigg(\vert DT_{\tau}(u_{1})\vert-(z_{2},DT_{\tau}(u_{1}))\bigg)\geq 0,$ $\displaystyle{\int_{\Omega}}\bigg(\vert DT_{\tau}(u_{2})\vert-(z_{1},DT_{\tau}(u_{2}))\bigg)\geq 0,$  employing \eqref{70} and \eqref{71} we have
\begin{equation*}
    \int_{\Omega}\frac{1}{(u_{1}u_{2})^{\gamma}}f\bigg(u^{\gamma}_{2}-u^{\gamma}_{1}\bigg)\bigg(T_{\tau}(u_{1})-T_{\tau}(u_{2})\bigg)dy= 0.
\end{equation*}
Therefore $T_{\tau}(u_{1})=T_{\tau}(u_{2}) $ a. e. in $\Omega $ for every $ \tau>0.$ In particular
\begin{equation*} 
u_{1}=u_{2} \hspace*{0.5cm}\mbox{ a. e. }\hspace*{0.2cm}in\hspace*{0.2cm} \Omega.
\end{equation*}
\begin{flushright}
    $\square$
\end{flushright}
\end{dem}
\section{Solutions for  $0\leq f\in L^{m}(\Omega)$}
 In this section, let us consider the following problem
\begin{equation}\label{94}
    \left\{
  \begin{array}{ll}
   -\Delta_{1}u=\frac{f}{u^{\gamma}}& \hbox{in $\Omega,$} \\
  \frac{Du}{\vert Du\vert}.\sigma+\frac{u}{\vert u\vert}+\vert u\vert^{s-1}u=0 & \hbox{on $\partial\Omega,$} \\
  \end{array}
\right.
\end{equation}
where $ f$ being nonnegative function in $ L^{m}(\Omega),$ $ 0<\gamma\leq 1$ and $s\geq 1.$
\par Now, we introduce the definition of solution to \eqref{94}.
\begin{defn}
Assume \( f \) is a nonnegative function in \( L^{m}(\Omega) \). We define \( u \) an element of  \( BV(\Omega) \) as a weak solution to problem \eqref{2}, if there exist $v\in L^{s+1}(\partial \Omega),$ $ \vartheta \in L^{\infty}(\partial \Omega) $ and $z\in \mathcal{DM}^{\infty}_{loc}(\Omega)$ with $\Vert z\Vert_{L^{\infty}(\Omega)^{m}}\leq 1$ such that
\begin{eqnarray}
   & &\frac{f}{u^{\gamma}}\in L^{1}_{loc}(\Omega), \nonumber \\ 
   & &\chi_{\{0<u\}} \in BV_{loc}(\Omega),\label{16} \\
   & &-\chi^{*}_{\{0<u\}}\textnormal{div }z =\frac{f}{u^{\gamma}} \mbox{in } \mathcal{D}^{\prime}(\Omega),\label{17}\\
   & & (z,DT_{\tau}(u))=\vert DT_{\tau}(u)\vert \mbox{ as measures in } \Omega \mbox{ for any } \tau>0,\label{18}\\
   & & T_{\tau}(u)+T_{\tau}(v)^{s+1}+T_{\tau}(u)[z,\sigma]=0\hspace*{0.4cm}  \mathcal{H}^{m-1}\mbox{ a. e. on } \partial \Omega \mbox{ for any } \tau>0, \label{20}\\
   & & \vartheta+v^{s}+[z,\sigma]=0 \hspace*{0.4cm}\mathcal{H}^{m-1} \mbox{ a. e. on } \partial \Omega. \label{45}
\end{eqnarray}\label{3}
\end{defn}
Our main existence result is stated as follows.
\begin{thm}\label{15}
Let $0\leq f\in L^{m}(\Omega),$  $ 0<\gamma\leq 1$ and $s\geq 1.$ Then there exists a  solution $ u$ to problem \eqref{4}.
\end{thm}
\begin{dem}
Let $ u_{p}$ be a solution to \eqref{55}, by Lemma \ref{5} there exist functions  $u$ belongs to $BV(\Omega),$ $\vartheta\in L^{\infty}(\partial \Omega)  $ and $ v\in L^{s+1}(\partial \Omega)$ satisfying, up to subsequences that $u_{p}$ converges strongly to $u $ in $ L^{1}(\Omega)$ as $ p\rightarrow 1^{+}$ and $u_{p}\rightharpoonup v $ weakly in $ L^{s+1}(\partial \Omega).$   Additionally, reasoning as in the proof of Lemma \ref{6} that there exists a vector field $z\in L^{\infty}(\Omega, \mathbb{R}^{m})$ with $\Vert z\Vert_{L^{\infty}(\Omega)^{m}}\leq 1$ such that  \eqref{18} hold. On the other hand choosing $\varphi\in C^{\infty}_{0}(\Omega)$ as a \textit{t. f.} in \eqref{55}, it yields
\begin{equation*}
    \int_{\Omega}\vert \nabla u_{p}\vert^{p-2}\nabla u_{p}\nabla \varphi dy=\int_{\Omega}\frac{f}{u^{\gamma}_{p}}\varphi dy.
\end{equation*}
From Fatou's Lemma we have
\begin{equation}\label{118}
\int_{\Omega}\frac{f}{u^{\gamma}}\varphi dy\leq \int_{\Omega}z\nabla \varphi dy=-\int_{\Omega}\varphi\textnormal{div }z,\hspace*{0.5cm}\forall\varphi\in C^{\infty}_{0}(\Omega).
\end{equation}
Then $z\in\mathcal{DM}_{loc}^{\infty}(\Omega)$ and $\frac{f}{u^{\gamma}}\in L_{loc}^{1}(\Omega).$ Moreover similarly to the proof of Lemma \ref{7}, it can be demonstrated that  \eqref{20} and that  \eqref{45} holds.
\par Now, we still have to prove \eqref{16} and \eqref{17},  we use $ G_{\delta}(u_{p})\varphi$ with $\varphi\in C^{\infty}_{0}(\Omega)$ as a \textit{t. f.} in \eqref{55}, with the support of Young's inequality, it appears that
\begin{equation}\label{30}
\int_{\Omega}\vert \nabla G_{\delta}(u_{p})\vert \varphi dy+\int_{\Omega}\vert \nabla u_{p}\vert^{p-2}\nabla u_{p} G_{\delta}(u_{p})\nabla \varphi dy\leq \frac{p-1}{p}\int_{\Omega} G_{\delta}^{\prime}(u_{p})\varphi dy+\int_{\Omega}\frac{f}{u^{\gamma}_{p}} G_{\delta}(u_{p})\varphi dy.
\end{equation}
We remark that \(  G_{\delta}(u_{p}) \) is bounded in \( BV(\Omega) \) regarding  \( p \), we apply lower semicontinuity, it is evident that
\begin{align}\label{138}
&\int_{\Omega}\varphi\vert DG_{\delta}(u)\vert\leq \lim_{p\rightarrow1^{+}}\int_{\Omega}\vert \nabla G_{\delta}(u_{p})\vert \varphi dy
\end{align}
Since \( G_{\delta}(u_{p})\overset{\ast}{\rightharpoonup} G_{\delta}(u) \) in \( L^{\infty}(\Omega) \) as \( p \rightarrow 1^{+} \), together with \eqref{200}, one observes that
\begin{align}\label{139}
&\lim_{p\rightarrow1^{+}}\int_{\Omega}\vert \nabla u_{p}\vert^{p-2}\nabla u_{p} G_{\delta}(u_{p})\nabla \varphi dy= \int_{\Omega}z G_{\delta}(u)\nabla \varphi dy.
\end{align}
With the help of Lebesgue's Theorem, it appears that
\begin{align}\label{140}
&\lim_{p\rightarrow1^{+}}\int_{\Omega}\frac{f}{u^{\gamma}_{p}} G_{\delta}(u_{p})\varphi dy=\int_{\Omega}\frac{f}{u^{\gamma}} G_{\delta}(u)\varphi dy
\end{align}
Passing to the limit as \( p \rightarrow 1^{+} \) in \eqref{30}, by \eqref{138}, \eqref{139} and \eqref{140}. Then
\begin{align}\label{114}
&\int_{\Omega}\varphi\vert DG_{\delta}(u)\vert +\int_{\Omega}z G_{\delta}(u)\nabla \varphi dy\leq\int_{\Omega}\frac{f}{u^{\gamma}} G_{\delta}(u)\varphi dy.
\end{align}
Since $ \frac{f}{u^{\gamma}}\in L^{1}_{loc}(\Omega),$ $ G_{\delta}(u)\leq1$  and  $z$ is a member of $\mathcal{DM}_{loc}^{\infty}(\Omega),$   we deduce that $ G_{\delta}(u)$ is  bounded in $BV_{loc}(\Omega)$ concerning $\delta.$
Taking the limit as $ \delta\rightarrow 0^{+}$ in \eqref{114}, by applying lower semicontinuity to the first term on the left-hand side and using Lebesgue's Theorem for the remaining terms, we find that
\begin{align}\label{116}
&\int_{\Omega}\varphi\vert D\chi_{\{0<u\}}\vert +\int_{\Omega}z\nabla \varphi \chi_{\{0<u\}}dy\leq\int_{\Omega}\frac{f}{u^{\gamma}}\varphi\chi_{\{0<u\}} dy.
\end{align}
Hence
$$\chi_{\{0<u\}} \in BV_{loc}(\Omega).$$
Notice that
\begin{equation}\label{115}
-\chi^{*}_{\{0<u\}}\textnormal{div }z=-\textnormal{div }(z\chi_{\{0<u\}})+(z,D\chi_{\{0<u\}}).
\end{equation}
Putting together $ \{u=0\}\subset \{f=0\}$ and \eqref{115}, then \eqref{116} becomes
\begin{equation}\label{117}
    -\int_{\Omega}\chi^{*}_{\{0<u\}}\varphi\textnormal{div }z\leq \int_{\Omega}\frac{f}{u^{\gamma}}\chi_{\{0<u\}}\varphi dy=\int_{\Omega}\frac{f}{u^{\gamma}}\varphi dy.
\end{equation}
We take $(\chi_{\{0<u\}}\ast\rho_{\epsilon})\varphi$ as a \textit{t. f.} in \eqref{118} with $0\leq\varphi\in C^{\infty}_{0}(\Omega)$ and $\rho_{\epsilon} $ is a mollifier  we arrive at
\begin{equation*}
\int_{\Omega}\frac{f}{u^{\gamma}}(\chi_{\{0<u\}}\ast\rho_{\epsilon})\varphi dy\leq -\int_{\Omega}(\chi_{\{0<u\}}\ast\rho_{\epsilon})\varphi \textnormal{div }z.
\end{equation*}
Passing to the limit $ \epsilon \rightarrow 0^{+}$ in previous inequality, on the left hand side we use Fatou's Lemma and apply Lebesgue Theorem on the right hand side. Then
\begin{align}\label{120}
&\chi_{\{0<u\}}\frac{f}{u^{\gamma}}=\frac{f}{u^{\gamma}}\leq-\chi^{\ast}_{\{0<u\}} \textnormal{div }z \hspace*{0.3cm}in \hspace*{0.3cm}\mathcal{D}^{\prime}(\Omega).
\end{align}
In accordance with \eqref{117} and \eqref{120} we obtain
\begin{align*}
&-\chi^{\ast}_{\{0<u\}} \textnormal{div }z=\frac{f}{u^{\gamma}}\hspace*{0.3cm}in \hspace*{0.3cm}\mathcal{D}^{\prime}(\Omega).
\end{align*}
\begin{flushright}
    $\square$
\end{flushright}
\end{dem}

\subsubsection*{Acknowledgements}
A heartfelt thank you to Professor Sergio Segura de Le\'{o}n for his invaluable assistance in revising and improving our paper.

\subsubsection*{Disclosure statement}
No potential conflict of interest was reported by the author(s).
\subsubsection*{Conflict of Interest}
The authors certify that they have no affiliations with or involvement in any organization or entity with any financial interest, or non financial interest (such as personal or professional relationships, knowledge or beliefs...).

\end{document}